\documentclass[aos,MSNbibl,seceqn,nameyear,dvips]{arximspdf}

\doi{10.1214/13-AOS1091} 
\volume{41}
\issue{2}
\pubyear{2013}
\firstpage{641}
\lastpage{669}

\makeatletter
\newcommand{\rrvert}{\vert}
\newcommand{\llvert}{\vert}

\renewcommand{\epsilon}{\varepsilon}

\newtheorem{theorem}{Theorem}[section]
\newtheorem{lem}{Lemma}[section]
\newtheorem{cor}{Corollary}[section]

\newproclaim{rmk}{Remark}[section]

\makeatother

\begin{document}
\begin{frontmatter}

\title{Moderate deviations for a nonparametric estimator of sample coverage}
\runtitle{MDP for a nonparametric estimator of sample coverage}

\begin{aug}
\author[A]{\fnms{Fuqing} \snm{Gao}\corref{}\thanksref{t1}\ead[label=e1]{fqgao@whu.edu.cn}}
\runauthor{F. Gao}
\affiliation{Wuhan University}
\address[A]{School of Mathematics and Statistics\\
Wuhan University\\
Wuhan 430072\\
P.R. China\\
\printead{e1}} 
\end{aug}

\thankstext{t1}{Supported by the National Natural Science
Foundation of China (10871153, 11171262).}

\received{\smonth{8} \syear{2011}}
\revised{\smonth{1} \syear{2013}}

%
\begin{abstract}
In this paper, we consider moderate deviations for Good's coverage
estimator. The moderate deviation principle and the self-normalized
moderate deviation principle for Good's coverage estimator are
established. The results are also applied to the hypothesis testing
problem and the confidence interval for the coverage.
\end{abstract}

%

\end{frontmatter}

\section{Introduction}\label{sec1}

Let $X_k(n)$ be the frequency of the $k$th species in
a random sample of size $n$ from a multinomial population with a
perhaps countably infinite number of species and let $ P_n$ be
probability measures under which the
$k$th species has probability $p_{kn}$ of being sampled, where $p_n
= (p_{kn}; k\geq1)$ with \mbox{$\sum_{k=1}^\infty p_{kn} = 1$}.
Let $Q_n$ and $F_j(n)$ denote the sum of the probabilities of the
unobserved species, and the total number of species represented $j$
times in the
sample, respectively, that is,
%
\begin{equation}
\label{def-Qn} Q_n =\sum_{k=1}^\infty
p_{kn}\delta_{k0}(n),\qquad F_j(n)=\sum
_{k=1}^\infty\delta_{kj}(n),
\end{equation}
where $ \delta_{kj}(n)= I_{\{X_k(n) = j\}}$. Then $1-Q_n$ is called the
sample coverage which is the sum of the probabilities of the observed
species. \citet{Goo53} proposed the estimator
%
\begin{equation}
\label{def-estimator}
\hat{Q}_n=\frac{F_{1}(n)}{n}
\end{equation}
for $Q_n $.

The Good estimator $ \hat{Q}_n$ has many applications such as
Shakespeare's general vocabulary and authorship of a poem
[\citet{EfrThi76}, \citet{ThiEfr87}], genom
[\citet{MaoLin02}], the probability of discovering new species\vadjust{\goodbreak} in
a population [\citet{GooTou56}, \citet{Cha81}], network
species and data confidentiality [\citet{Zha05}].
\citet{LlaGouRee11} considered the problem of predicting $Q_n$.
They studied prediction and prediction intervals, and gave a real-data
example.

On the theoretical aspects, many authors studied the asymptotic
properties [cf. Esty (\citeyear{Est82,Est83}),
\citet{OrlSanZha03}, and \citet{ZhaZha09} and references
therein]. \citet{Est83} proved the following asymptotic normality:
%
\begin{equation}
\label{clt-thm-esty} \lim_{n\to\infty}P_n \biggl(
\frac{n(\hat{Q}_{n}-Q_n)}{\sqrt {b(n)}}\leq x \biggr)=\int_{-\infty}^x
\frac{1}{\sqrt{2\pi}}e^{-{u^2}/{2}}\,du,\qquad x\in\mathbb R,
\end{equation}
under the condition
%
\begin{equation}
\label{esty-condition}\qquad \lim_{n\to\infty} \frac{E_n(F_1(n))}{n}=c_1
\in(0,1)\quad\mbox{and}\quad\lim_{n\to\infty} \frac{E_n(F_2(n))}{n}=c_2
\in[0,\infty),
\end{equation}
where
%
\begin{equation}
\label{constant-bn} b(n)= E_n\bigl(F_1(n)\bigr)
\bigl(1-E_n\bigl(F_1(n)\bigr)/n\bigr)+2E_n
\bigl(F_2(n)\bigr).
\end{equation}
Recently, \citet{ZhaZha09} found a necessary and sufficient condition
for the asymptotic normality (\ref{clt-thm-esty}) under the condition
%
\begin{equation}
\label{moment-condition-1} \limsup_{n\to\infty} \frac{E_n(F_1(n))}{n}<1,
\end{equation}
that is, under condition (\ref{moment-condition-1}), (\ref
{clt-thm-esty}) holds if and only if both
%
\begin{equation}
\label{moment-condition-2} \lim_{n\to\infty} \bigl( E_n
\bigl(F_1(n)\bigr) + E_n\bigl(F_2(n)\bigr)
\bigr)=\infty
\end{equation}
and for any $\varepsilon>0$,
%
\begin{equation}
\label{clt-Lindeberg-condition} \lim_{n\to\infty} \frac{1}{ s^2_{n}}
\sum_{k=1}^\infty(n p_{kn})^2e^{-n
p_{kn}}
I_{\{n
p_{kn}>\varepsilon s_{n}\}}=0,
\end{equation}
where for any $\lambda>0$,
%
\begin{equation}
\label{var-def} s_{\lambda n}^2=\sum
_{k=1}^\infty\bigl(\lambda p_{kn}e^{-\lambda p_{kn}}
+(\lambda p_{kn})^2e^{-\lambda p_{kn}}\bigr)\quad \mbox{and}\quad
s_n=s_{nn}.
\end{equation}

In this paper, we consider the moderate deviation problem for the Good
estimator. It is known that the moderate deviation principle is a basic
problem. It provides us with rates of convergence and a useful method
for constructing asymptotic
confidence intervals. The moderate deviations can be applied to the
following nonparameter hypothesis testing problem:
\[
H_0\dvtx  P_n=P_n^{(0)} \quad\mbox{and}\quad
H_1\dvtx  P_n=P_n^{(1)},
\]
where $P_n^{(0)}$ and $P_n^{(1)}$ are two probability measures under
which the
$k$th species has, respectively, probability $p_{kn}^{(0)}$ and
$p_{kn}^{(1)}$ of being sampled, where $p_n^{(i)}
= (p_{kn}^{(i)}; k\geq1)$ with $\sum_{k=1}^\infty p_{kn}^{(i)} = 1$,
$i=0,1$. We can define a rejection region of\vadjust{\goodbreak} the hypothesis testing by
the moderate deviation principle such that the probabilities of type I
and type II errors tend to $0$ with an exponential speed. The
asymptotic normality provides $\sqrt{b(n)}$ as the asymptotic variance
and approximate confidence statements, but it does not prove that the
probabilities of type I and type II errors tend to $0$ with an
exponential speed. The moderate deviations can be applied to a
hypothesis testing problem for the expected coverage of the sample.

\citet{GaoZha11} have established a general delta method on the
moderate deviations for estimators. But the method cannot be applied to
the Good estimator. In order to study the moderate deviation problem
for the Good estimator, we need refined asymptotic analysis techniques
and tail probability estimates. The exponential moments inequalities,
the truncation method, asymptotic analysis techniques and the Poisson
approximation in \citet{ZhaZha09} play important roles. Our main
results are a moderate deviation principle and a self-normalized
moderate deviation principle for the Good estimator.

The rest of this paper is organized as follows. The main results are
stated in Section~\ref{sec2}. Some examples and applications to the hypothesis
testing problem and the confidence interval are also given in Section
\ref{sec2}. The proofs of the main results are given in Section~\ref{sec3}. Some basic
concepts for large deviations and the proofs of several technique
lemmas are given in the \hyperref[app]{Appendix}.

\section{Main results and their applications}\label{sec2}

In this section, we state the main results and give some examples and
applications.

\subsection{Main results}\label{sec2.1}

Let $a(t)$, $t\ge0$, be a function taking values in $[1,+\infty)$ such that
%
\begin{equation}
\label{constant-an-1} \lim_{t\to\infty} \frac{a(t)}{\sqrt{t}} =
\infty,\qquad
\lim_{t\to\infty}\frac{ a(t)}{t} =0.
\end{equation}

We introduce the following Lindeberg-type condition: for any positive
sequence $\{\lambda_n,n\geq1\}$ with $\lambda_n/n\to1$ and any
$\varepsilon>0$,
%
\begin{equation}
\label{mdp-Lindeberg-condition} \lim_{n\to\infty} \frac{1}{s_{n}^2}
\sum_{k=1}^\infty(\lambda _np_{kn})^2e^{- \lambda_n p_{kn}}
I_{\{\lambda_np_{kn}>\varepsilon
s_{n}^2/a(s_{n}^2)\}}=0.
\end{equation}

\begin{rmk}\label{rmk-2-1}
For any $L\ge1$,
\begin{eqnarray*}
&&\sum_{k=1}^\infty(\lambda_n
p_{kn})^2e^{-\lambda_n p_{kn}} I_{\{
\lambda_n
p_{kn}> L \}}
\\
&&\qquad\leq \sum_{j=0}^\infty L 2^{j+1}
\exp \bigl\{- L 2^j \bigr\} \sum_{k=1}^\infty
\lambda_n p_{kn} I_{\{ L 2^j \leq\lambda_n p_{kn}< L 2^{j+1} \}}\\
&&\qquad\leq8
\lambda_n L \exp \{- L \}.
\end{eqnarray*}
In particular, take $L=\frac{\varepsilon s_{n}^2}{a(s_{n}^2)}$. If
$\lim_{n\to\infty} \frac{s_{n}^2}{a(s_{n}^2)\log(\lambda
_n/s_n^2)}=\infty$,
then (\ref{mdp-Lindeberg-condition}) holds.\vadjust{\goodbreak}
\end{rmk}

\begin{theorem}[(Moderate deviation principle)]\label{main-thm-mdp}
Suppose that the conditions (\ref{moment-condition-1}), (\ref
{moment-condition-2}) and (\ref{mdp-Lindeberg-condition}) hold. Then
$ \{\frac{n(\hat{Q}_{n}-Q_n)}{a(b(n))},n\geq1 \}$
satisfies a large deviation principle with speed
$\frac{a^2(b(n))}{b( n)}$ and with rate function $I(x)=\frac{x^2}{2}$.
In particular, for any $r>0$,
\[
\lim_{n\to\infty}\frac{b( n)}{a^2(b(n))}\log P_n \biggl(
\frac
{|n(\hat
{Q}_{n}-Q_n)|}{a(b(n))} \geq r \biggr)= -\frac{r^2}{2}.\vspace*{-2pt}
\]
\end{theorem}

\begin{theorem}[(Self-normalized moderate deviation
principle)]\label{main-thm-mdp-self}
Suppose that conditions (\ref{moment-condition-1}), (\ref
{moment-condition-2}) and (\ref{mdp-Lindeberg-condition}) hold. Then
\[
\biggl\{\frac{\sqrt{b(n)}n(\hat{Q}_{n}-Q_n)}{a(b(n)) \sqrt {F_1(n)(1-F_1(n)/n)+2F_2(n) }},n\geq1 \biggr\}
\]
satisfies a large deviation principle with speed
$\frac{a^2(b(n))}{b( n)}$ and with rate function $I(x)=\frac{x^2}{2}$.\vspace*{-2pt}
\end{theorem}

\begin{rmk}
Let $t_n$, $n\ge1$ be a sequence of positive numbers such that
%
\begin{equation}
\label{constant-tn-1} t_n\uparrow\infty\quad\mbox{and}\quad
\frac{t_n}{\sqrt {b(n)}}\downarrow0.
\end{equation}
Then Theorems~\ref{main-thm-mdp} and~\ref{main-thm-mdp-self}
give the following estimates which are much easier to understand and apply:
\[
P_n \biggl(\pm\frac{ n(\hat{Q}_{n}-Q_n) }{ \sqrt{b(n)} } \geq t_n \biggr)=
\exp \biggl\{-\bigl(1+o(1)\bigr)\frac{t_n^2}{2} \biggr\}
\]
and
\[
P_n \biggl(\pm\frac{ n(\hat{Q}_{n}-Q_n) }{ \sqrt {F_1(n)(1-F_1(n)/n)+2F_2(n) }} \geq t_n \biggr)=
\exp \biggl\{ -\bigl(1+o(1)\bigr)\frac{t_n^2}{2} \biggr\}.
\]
\end{rmk}

Set $u_n=E_n(Q_n)=\sum_{k=1}^\infty p_{kn}(1-p_{kn})^{n}$. Then $1-u_n$
is called the expected coverage\vspace*{1pt} of the sample in the literature. By
Theorems~\ref{main-thm-mdp} and~\ref{main-thm-mdp-self}, and
Lemma~\ref{main-thm-self-lem}, $ \hat{Q}_{n} $ as an estimator of $u_n$
also satisfies moderate deviation principles.\vspace*{-2pt}

\begin{cor}\label{main-thm-mdp-expected}
Suppose that conditions (\ref{moment-condition-1}), (\ref
{moment-condition-2}) and (\ref{mdp-Lindeberg-condition}) hold. Then
$ \{\frac{n(\hat{Q}_{n}-u_n)}{a(b(n))},n\geq1 \}$ and $
\{\frac{\sqrt{b(n)}n(\hat{Q}_{n}-u_n)}{a(b(n)) \sqrt {F_1(n)(1-F_1(n)/n)+2F_2(n) }},n\geq1 \}
$
satisfy the large deviation principle with speed
$\frac{a^2(b(n))}{b( n)}$ and with rate function \mbox{$I(x)=\frac{x^2}{2}$}.\vspace*{-2pt}
\end{cor}

\begin{rmk}
\citet{LlaGouRee11} considered the problem of predicting
$Q_n$, and obtained conditionally unbiased predictors and exact
prediction intervals based on a Poissonization argument. The moderate
deviations for the predictors are also interesting problems.\vadjust{\goodbreak}
\end{rmk}

\subsection{Application to hypothesis testing and confidence interval}\label{sec2.2}

In this subsection, we apply the moderate deviations to hypothesis
testing problems and confidence interval.
Let $Q_n$ be the unknown total probability unobserved species, and
let $\hat{Q}_n$ be the estimator defined by (\ref{def-estimator}).

First, let us consider a nonparametric hypothesis
testing problem. Let $P_n^{(0)}$ and $P_n^{(1)}$ be two probability
measures under which the
$k$th species has, respectively, probability $p_{kn}^{(0)}$ and
$p_{kn}^{(1)}$ of being sampled, where $p_n^{(i)}
= (p_{kn}^{(i)};\break k\geq1)$ with $\sum_{k=1}^\infty p_{kn}^{(i)} = 1$,
$i=0,1$. Denote by
\[
u_n^{(i)}:= \sum_{k=1}^\infty
p_{kn}^{(i)} \bigl(1-p_{kn}^{(i)}
\bigr)^{n},\qquad i=0,1,
\]
and
\[
b^{(i)}(n):= E_n^{(i)}\bigl(F_1(n)
\bigr) \bigl(1-E_n^{(i)}\bigl(F_1(n)\bigr)/n
\bigr)+2E_n^{(i)}\bigl(F_2(n)\bigr),\qquad i=0,1.
\]
Suppose that the conditions (\ref{moment-condition-1}), (\ref
{moment-condition-2}) and (\ref{mdp-Lindeberg-condition}) hold for
$P_n^{(i)}$, $i=0,1$, and that
\[
\liminf_{n\to\infty} \bigl| u_n^{(0)}-
u_n^{(1)}\bigr|\not=0.
\]
Consider the nonparameter hypothesis testing
\[
H_0\dvtx  P_n=P_n^{(0)} \quad\mbox{and}\quad
H_1\dvtx  P_n=P_n^{(1)}.
\]
We take the
statistic
$
T_n:= \hat{Q}_n-u_n^{(0)}
$
as test statistic. Suppose
that the rejection region for testing the null hypothesis $H_0$
against $H_1$ is \mbox{$\{\frac{n}{a(b^{(0)}(n))}|T_n|\geq c\}$}, where $c$
is a
positive constant. The probability $\alpha_n$ of type I error
and the probability $\beta_n$ of type II error are
\[
\alpha_n=P_n^{(0)} \biggl(
\frac{n}{a(b^{(0)}(n))}|T_n|\geq c \biggr),\qquad \beta_n=P_n^{(1)}
\biggl(\frac{n}{a(b^{(0)}(n))}|T_n|< c \biggr),
\]
respectively. It follows
\[
\beta_n \leq P_n^{(1)} \biggl(
\frac{n}{a(b^{(1)}(n))}\bigl\llvert \hat {Q}_n-u_n^{(1)}
\bigr\rrvert \geq \biggl( \bigl| u_n^{(0)}- u_n^{(1)}\bigr|
-\frac{a(b^{(0)}(n))c}{n} \biggr)\frac
{n}{a(b^{(1)}(n))} \biggr).
\]
Therefore, Corollary~\ref{main-thm-mdp-expected} implies that
\[
\lim_{n\to\infty}\frac{b^{(0)}(n)}{a^2(b^{(0)}(n))}\log\alpha_n=-
\frac{c^2}{2},\qquad \lim_{n\to\infty}\frac
{b^{(1)}(n)}{a^2(b^{(1)}(n))}\log
\beta_n=- \infty.
\]

The above result tells us that if the rejection region for the test
is\break $\{\frac{n}{a(b^{(0)}(n))}|T_n|\geq c\}$, then the probability of type
I error tends to $0$ with exponential decay speed
$\exp\{- c^2 a^2(b^{(0)}(n))/(2b^{(0)}(n))\}$, and the probability of
type II error tends to $0$ with exponential decay speed $\exp\{- r
a^2(b^{(1)}(n))/\break b^{(1)}(n) \}$ for all $r>0$. But the asymptotic
normality does not prove that the probabilities of type I and type II
errors tend to $0$ with an exponential speed.\vadjust{\goodbreak}

We also consider a hypothesis testing problem for the expected coverage
of the sample. We denote by $P_n^{u_n}$ the probability measures under
which the expected coverage of the sample is $1-u_n$ and set
\[
b_{u_n}(n):= E_n^{u_n}\bigl(F_1(n)
\bigr) \bigl(1-E_n^{u_n}\bigl(F_1(n)\bigr)/n
\bigr)+2E_n^{u_n}\bigl(F_2(n)\bigr).
\]
Suppose that the conditions (\ref{moment-condition-1}), (\ref
{moment-condition-2}) and (\ref{mdp-Lindeberg-condition}) hold for
$P_n^{u_n}$ for each $u_n>0$. Let $0< u_n^{(0)}\leq u_n^{(1)}$ be two
real numbers preassigned. Consider the hypothesis testing
\[
H_0\dvtx  u_n\leq u_n^{(0)}
\quad\mbox{and}\quad H_1\dvtx  u_n> u_n^{(1)}.
\]
We also take the rejection region $D_
n:=\{\frac{n}{a(b_{u_n^0}(n))}(\hat{Q}_n-u_n^{(0)})\geq c\}$, where
$c$ is a
positive constant. When $u_n\leq u_n^{(0)}$,
\[
\log P_n^{u_n} (D_n )\leq\log
P_n^{u_n} \biggl(\frac
{n}{a(b_{u_n^0}(n))}(\hat{Q}_n-u_n)
\geq c \biggr)\approx- \frac{ c^2
a^2(b_{u_n^0}(n))}{2b_{u_n}(n)}
\]
and when $u_n> u_n^{(1)}$,
\[
\frac{b_{u_n}(n)}{a^2(b_{u_n}(n))}\log P_n^{u_n} \bigl(D_n^c
\bigr) \rightarrow- \infty.
\]

Next, we apply the moderate estimates to confidence intervals. For
given confidence level
$1-\alpha$, set $c_\alpha=\sqrt{- \frac{b(n)}{a^2(b(n))}\log
\alpha}$. Then by Theorem~\ref{main-thm-mdp}, the $1-\alpha$ confidence
interval for $Q_n$ is
$ (\hat{Q}_n- \frac{a(b(n))}{n}c_\alpha,\hat{Q}_n+ \frac
{a(b(n))}{n}c_\alpha)
$, that is,
\[
\biggl(\hat{Q}_n-\frac{1}{n} \sqrt{- b(n)\log \alpha},
\hat{Q}_n+ \frac{1}{n} \sqrt{- b(n)\log \alpha} \biggr).
\]
But the confidence interval contains unknown $b(n)$. We use Theorem
\ref
{main-thm-mdp-self} to obtain another confidence interval with
confidence level
$1-\alpha$ for $Q_n$ which does not contain unknown $b(n)$,
\begin{eqnarray*}
&& \biggl(\hat{Q}_n-\frac{\sqrt{- (F_1(n)(1-F_1(n)/n)+2F_2(n))\log
\alpha}}{n},
\\
&&\qquad\hspace*{0pt} \hat{Q}_n+ \frac{\sqrt{-
(F_1(n)(1-F_1(n)/n)+2F_2(n))\log
\alpha}}{n} \biggr).
\end{eqnarray*}

\subsection{Examples}\label{sec2.3} Let us check that some examples
in \citet{ZhaZha09} also satisfy moderate deviation principles if
$a(n)=n^{\gamma}$, where $\gamma\in(1/2,1)$. For a given decreasing
density function $p_n(x)$ on $[0,\infty)$. Define
$p_{in}=z_np_n(i)$, where $z_n=(\sum_{i=1}^\infty p_{in})^{-1}$.
Two concrete examples are as follows:

Let $p_n(x)=p(x)=a/(x+1)^b$, where $a>0$ and $b>1$. By Example 1 in
\citet{ZhaZha09},
$E_n(F_1(n))\asymp n^{1/b}$ and $\log s_n^2\asymp\log n$,
where
\[
c_n\asymp b_n \mbox{ means } 0<\liminf
_{n\to\infty}\frac
{c_n}{b_n}\leq \limsup_{n\to\infty}
\frac{c_n}{b_n}<\infty.\vadjust{\goodbreak}
\]
Thus (\ref{moment-condition-1}) and (\ref{moment-condition-2}) hold. By
Remark~\ref{rmk-2-1}, (\ref{mdp-Lindeberg-condition}) also holds.
Therefore, Theorems~\ref{main-thm-mdp} and~\ref{main-thm-mdp-self} hold.

Let $p_n(x)=p(x)=r_n^{-1}e^{-x/r_n}$, where $r_n/n\leq c $ for some
constant $c<\infty$. Then by Example 2 in \citet{ZhaZha09},
$\limsup_{n\to\infty}\frac{E_n(F_1(n))}{n}= \limsup_{n\to\infty
}\int_0^1 e^{-n y/r_n}\,dy\leq\int_0^1 e^{- y/c}\,dy<1$
and $s_{\lambda_nn}^2\asymp r_n\int_0^{\lambda
_n/r_n}(1+t)e^{-t}\,dt\asymp r_n$ when $\lambda_n/n\to1$.
Thus, (\ref{mdp-Lindeberg-condition}) is equivalent to
\[
o(1)=\frac{1}{r_n}\int_{n p_n(x)\geq\varepsilon r_n/a(r_n)} \bigl(
\lambda_n p_n(x)\bigr)^2 e^{-\lambda_np_n(x)}\,dx=
\int_{\varepsilon r_n/a(r_n) \leq
t\leq\lambda_n/r_n} t e^{-t} \,dt,
\]
which holds if and only if $r_n\to\infty$.
Therefore, (\ref{moment-condition-1}), (\ref{moment-condition-2}) and
(\ref{mdp-Lindeberg-condition}) hold if and only if $r_n\to\infty$.

\section{Proofs of main results}\label{sec3}

In this section we give proofs of the main results.
Let us explain the idea of the proof of Theorem~\ref{main-thm-mdp}.
First, we divide the proof into two cases: case I and case II,
according to the limit $\lim_{n\to\infty}{E_n(F_1(n))}/{n} \in(0,1)$
and $0$.
For case I, by the truncation method and the exponential equivalent
method, we simplify our problems to the case which $\{np_{nk},k\geq
1,n\geq1\}$ is uniformly bounded. For case II, by the Poisson
approximation and the exponential equivalent method, we simplify our
problems to the case of independent sums satisfying an analogous
Lindeberg condition. For the two cases simplified, we establish
moderate deviation principles by the method of the Laplace asymptotic
integral (Lemmas~\ref{laplace-int-lem} and~\ref{mdp-poisson-approx}).
The exponential moment estimate (Lemma~\ref{exp-moment-inq-lem-1})
plays an important role in the proofs of some exponential equivalence
(Lemmas~\ref{exp-moment-estimates-lem-1} and~\ref{exp-eqv-lemma-1}).
The main technique in the estimate of the Laplace asymptotic integral
Lemma~\ref{laplace-int-lem} is asymptotic analysis. In particular, we
emphasis a transformation defined below (\ref{Laplace-int-lem-eq-4})
which plays a crucial role in the proof of Lemma~\ref{laplace-int-lem}.

We can assume that the population is
sampled sequentially, so that ${\mathbf X}(m)-{\mathbf X}(m-1)$,
$m\geq1$, are i.i.d. $\operatorname{multinomial}(1, p_n)$ under $P_n$, where
${\mathbf X}(n) = (X_k(n)$; $k\geq1)$ can be viewed
as a multinomial $(n; p_n)$ vector under $P_n$, that is, for all
integers $m\geq1$,
\[
P_n \bigl(X_k(n) = x_k; k=1, \ldots,m
\bigr)= \frac{n!(1-\sum_{k=1}^m
p_{kn})^{n-x_1-\cdots-x_m}\prod_{k=1}^mp_{kn}^{x_k}}{(n-x_1-\cdots
-x_m)!x_1!\cdots
x_m!}.
\]

It is obvious that $ E_n(F_1(n)) /n\leq1 $. Since for any $1\leq L< n$,
\begin{eqnarray*}
\frac{2E_n(F_2(n))}{n-1}&\leq& L\sum_{np_{kn}\leq
L}p_{kn}(1-p_{kn})^{n-2}+
\sup_{np\geq L}np(1-p)^{n-2}
\\
&\leq& \frac{L}{1-L/n} \frac{E_n(F_1(n))}{n}+Le^{-L} \biggl(1-
\frac
{L}{n} \biggr)^{-2},
\end{eqnarray*}
we have that
%
\begin{equation}
\label{mdp-condition-4} \limsup_{n\to\infty} \frac{2E_n(F_2(n))}{n}\leq
\limsup_{n\to
\infty} \frac{E_n(F_1(n))}{n}+e^{-1} \leq2;
\end{equation}
and if $ \limsup_{n\to\infty} \frac{E_n(F_1(n))}{n}=0$, then $
\limsup_{n\to\infty} \frac{E_n(F_2(n))}{n}=0$.
Without loss of generality, we can assume that
%
\begin{equation}
\label{mdp-condition-5}\quad \lim_{n\to\infty} \frac{E_n(F_1(n))}{n}=c_1
\in[0,1) \quad\mbox{and}\quad \lim_{n\to\infty} \frac{ E_n(F_2(n))}{n}=c_2
\in[0,1].
\end{equation}
Otherwise, we consider subsequence.
The proof of Theorem~\ref{main-thm-mdp} will be divided into two cases,
\[
\mbox{case I:}\quad  c_1\in(0,1);\qquad \mbox{case II:}\quad
c_1=0.
\]

Now let us introduce
the structrue of the proofs of main results. In Section~\ref{sec3.1}, we give
several moment estimates and exponential moment inequalities which are
basic for studying the moderate deviations for the Good estimator. A
truncation method and some related estimates are also presented in the
subsection. The proofs of cases I and II of Theorem \ref
{main-thm-mdp} are given, respectively, in Sections~\ref{sec3.2} and
\ref{sec3.3}. In Section~\ref{sec3.4}, we prove Theorem~\ref{main-thm-mdp-self}. The
proofs of several technique lemmas are postponed to the
\hyperref[app]{Appendix}.

\subsection{Several moment estimates and inequalities}\label{sec3.1}

For any $L\geq1$ and $\varrho>0$, set
\[
M_n^L= \{k\geq1; np_{kn}\leq L \},\qquad
M_n^{Lc}=\{k\geq 1;np_{kn}> L\}
\]
and
\begin{eqnarray*}
M_{n\varrho}&=&\bigl\{k\geq1;np_{kn}\leq\varrho b(n)/a(b(n)\bigr\},\\
M_{n\varrho}^{c}&=&\bigl\{k\geq1;np_{kn}>\varrho
b(n)/a(b(n)\bigr\}.
\end{eqnarray*}

\begin{lem}\label{truncation-lem-L}
If $c_1\in(0,1)$, then for any positive sequence $\{\lambda_n,n\geq
1\}
$ with $\lambda_n/n\to1$,
%
\begin{equation}
\label{truncation-lem-L-eq-1} \lim_{L\to\infty}\limsup
_{n\to\infty}\frac{1}{n} \sum_{k\in M_n^{Lc}}
\bigl(\lambda_n p_{kn}+(\lambda_n
p_{kn})^2 \bigr)e^{-\lambda_n
p_{kn}} =0.
\end{equation}
In particular, condition (\ref{mdp-Lindeberg-condition}) is
valid.
\end{lem}

\begin{pf}
Similarly to Remark~\ref{rmk-2-1}, for any $L\ge1$,
\begin{eqnarray*}
\sum_{k\in M_n^{Lc}} \lambda_n p_{kn}
e^{-\lambda_n p_{kn}} &\leq& \lambda_n e^{-L}/
\bigl(1-e^{-L}\bigr),
\\
\sum_{k\in M_n^{Lc}} (\lambda_n
p_{kn})^2e^{-\lambda_n p_{kn}} &\leq&8L \lambda_n\exp
\{- L \}.
\end{eqnarray*}
Therefore, (\ref{truncation-lem-L-eq-1}) holds.
\end{pf}

\begin{rmk}
From Lemma 1 in \citet{ZhaZha09}, under conditions (\ref
{moment-condition-1}) and (\ref{moment-condition-2}),
\[
\frac{E_n(F_1(n))+ 2 E_n(F_2(n))}{s_n^2}\to1,\qquad b(n) \asymp s_n^2,
\]
and if $c_1\in(0,1)$, then $
\lim_{n\to\infty}\frac{s_{n}^2}{n}=c_1+2c_2>0$.
\end{rmk}

\begin{lem}\label{sn-comparison-lem-1}
Assume that (\ref{mdp-Lindeberg-condition}) holds. If $0<\lambda
_n\leq
n$ and
\[
\limsup_{n\to\infty}\frac{n-\lambda
_n}{na(b(n))/b(n)}<\infty,
\]
then
%
\begin{equation}
\label{sn-comparison-lem-1-eq-1} s_{\lambda_n n}^2=\bigl(1+o(1)
\bigr)s_n^2.
\end{equation}
\end{lem}
\begin{pf} Set $r:=\limsup_{n\to\infty}\frac{n-\lambda
_n}{na(b(n))/b(n)}$. Then for any $\epsilon>0$, for $n$ large enough,
\begin{eqnarray*}
s_{\lambda_n n}^2 &\leq& e^\epsilon\sum
_{k=1}^\infty \bigl(n p_{kn}+(n
p_{kn})^2 \bigr)e^{-n p_{kn}}
\\
&&{} +\sum_{k=1}^\infty \bigl(
\lambda_n p_{kn}+(\lambda_n
p_{kn})^2 \bigr)e^{-n p_{kn}}I_{\{np_{kn}>\epsilon b(n)/(2r a(b(n)))\}}.
\end{eqnarray*}
Therefore, by (\ref{mdp-Lindeberg-condition}), the above inequality
implies that
$\limsup_{n\to\infty}\frac{ s_{\lambda_n n}^2}{s_n^2}\leq
e^\epsilon\to
1$ as $ \epsilon\to0$.
On the other hand, it is clear that for any $\epsilon>0$, when $n$ is
large enough,
\[
s_{\lambda_n n}^2\geq\sum_{k=1}^\infty
\bigl(\lambda_n p_{kn}+(\lambda_n
p_{kn})^2 \bigr)e^{-n p_{kn}} \geq(1-
\epsilon)^2 s_n^2,
\]
which yields that $\liminf_{n\to\infty}\frac{ s_{\lambda_n
n}^2}{s_n^2}\geq1$. Thus (\ref{sn-comparison-lem-1-eq-1}) is valid.
\end{pf}

\begin{lem}\label{truncation-lem-rho}
For any $\varrho>0$,
%
\begin{equation}
\label{truncation-lem-rho-eq-1}\quad \lim_{n\to\infty} \frac{1}{b(n)}
\sum_{k\in M_{n\varrho}} \biggl\llvert E_n\bigl(
\delta_{kj}(n)\bigr)- \frac
{1}{j!}(np_{kn})^j
e^{-np_{kn}}\biggr\rrvert=0,\qquad j=1,2.
\end{equation}
\end{lem}

\begin{pf}
Since $(1-p_{kn})^{n-j}= e^{-np_{kn}}(1+O(b(n)/a^2(b(n))))$ holds
uniformly on $M_{n\varrho}$ for $j=1,2$, we obtain that
\begin{eqnarray*}
\hspace*{-3pt}&&
\frac{1}{b(n)}\sum_{k\in M_{n\varrho}} \biggl\llvert
\frac{n!}{(n-j)!} p_{kn}^j (1-p_{kn})^{n-j}-(np_{kn})^j
e^{-np_{kn}}\biggr\rrvert
\\
\hspace*{-3pt}&&\qquad=\frac{1}{b(n)}\sum_{k\in M_{n\varrho}}(np_{kn})^j
e^{-np_{kn}} \biggl\llvert \frac{n!}{(n-j)!n^j}-\bigl(1+O\bigl(b(n)/a^2
\bigl(b(n)\bigr)\bigr)\bigr) \biggr\rrvert =o(1).
\end{eqnarray*}
That is, (\ref{truncation-lem-rho-eq-1}) holds.
\end{pf}

In order to obtain the exponential moment inequalities, we need some
concepts of negative dependence; cf. \citet{JoaPro83},
\citet{DubRan98}. Let $\eta_1,\eta_2,\ldots$ be real random
variables. $\eta_1,\eta_2,\ldots$ are said to be negatively associated
if for every two disjoint index finite sets $\Lambda_1,\Lambda
_2\subset
\{1,2,\ldots\}$,
\[
E \bigl(f(\eta_k,k\in\Lambda_1)g(\eta_k,k
\in\Lambda_2) \bigr)\leq E \bigl(f(\eta_k,k\in
\Lambda_1) \bigr)E \bigl(g(\eta_k,k\in \Lambda
_2) \bigr)
\]
for all nonnegative functions $f\dvtx \mathbb R^{\Lambda_1}\to\mathbb R$
and $g\dvtx \mathbb R^{\Lambda_2}\to\mathbb R$ that are both nondecreasing
or both nonincreasing.

\begin{lem}\label{negative-dep-lem}
$\{X_k(n),k\geq1\}$ is a sequences of negatively associated random
variables, and for each $0\leq j\leq n$ $\{\delta_{k0}(n) +\delta
_{k1}(n)+\cdots+\delta_{kj}(n),k\geq1\}$ is also negatively associated.
\end{lem}

\begin{pf}
Let $\delta_{k}^m$ denote the frequency of the $k$th species in the
$m$th sampling, that is,
\[
\delta_{k}^m=I_{\{X_k(m)-X_k(m-1)=1\}}.
\]
Then $\delta_{k}^m,k\geq1$ are zero-one random variables such that
$\sum_{k=1}^\infty\delta_{k}^m=1 $. By Lemma~8 in
\citet{DubRan98}, $\delta_{k}^m,k\geq1$, are negative associated.
Since $\{\delta_{k}^m,k\geq1\}$, $m= 1,\ldots,n$, are i.i.d. under
$P_n$, $ \delta_{k}^m,k\geq1$, $m= 1,\ldots,n $, are negative
associated. Noting that $X_k(n)=\sum_{m=1}^n \delta_{k}^m$ and
\[
\delta_{k0}(n) +\delta_{k1}(n)+\cdots+\delta_{kj}(n)=
\psi\bigl(X_k(n)\bigr),
\]
where $\psi(x)=I_{(-\infty, j]}(x)$ is a decreasing function, we obtain
that $\{X_k(n),\break k\geq1\}$ and $\{\delta_{k0}(n) +\delta
_{k1}(n)+\cdots
+\delta_{kj}(n),k\geq1\}$ are two sequences of negatively associated
random variables.
\end{pf}

\begin{lem}\label{exp-moment-inq-lem-1} Let $M$ be a subset of the set
$\mathbb N$ of positive integers. Then
for any $r\in\mathbb R$,
%
\begin{equation}
\label{exp-moment-inq-lem-eq-1} E_n \biggl(\exp \biggl\{r\sum
_{k\in M} p_{kn}\delta_{k0}(n)
\biggr\} \biggr) \leq\prod_{k\in M} \bigl(
\bigl(e^{rp_{kn}} -1 \bigr) (1-p_{kn})^n+1 \bigr)
\end{equation}
and for any $j\geq1$,
%
\begin{eqnarray}
\label{exp-moment-inq-lem-eq-2}
&&
E_n \biggl(\exp \biggl\{r\sum
_{k\in M} \bigl(\delta_{k0}(n) +\delta
_{k1}(n)+\cdots+\delta_{kj}(n) \bigr) \biggr\} \biggr)
\nonumber\\[-8pt]\\[-8pt]
&&\qquad\leq\prod_{k\in M} \Biggl( \bigl(e^r-1
\bigr)\sum_{l=0}^j \frac
{n!}{(n-l)!l!}p_{kn}^l
(1-p_{kn})^{n-l}+1 \Biggr).
\nonumber
\end{eqnarray}
\end{lem}

\begin{pf}
For any $r\in\mathbb R$ given, set $\psi_{k}(x)=e^{r p_{kn} x}$,
$x\in
\mathbb R$. Then, when \mbox{$r\geq0$}, all $\psi_{k}$, $k\geq1$ are
nonnegative and increasing; when $r< 0$, all $\psi_{k}$, $k\geq1$ are
nonnegative and decreasing. Therefore, by Lemma~\ref{negative-dep-lem},
\begin{eqnarray*}
E_n \biggl(\exp \biggl\{r\sum_{k\in M}
p_{kn}\delta_{k0}(n) \biggr\} \biggr) &\leq& \prod
_{k\in M}E_n \bigl(\exp \bigl\{r p_{kn}
\delta _{k0}(n) \bigr\} \bigr)
\\
&\leq&\prod_{k\in M} \bigl( \bigl(e^{rp_{kn}} -1
\bigr) (1-p_{kn})^n+1 \bigr).
\end{eqnarray*}
Similarly, we can obtain (\ref{exp-moment-inq-lem-eq-2}).
\end{pf}

As applications of Lemma~\ref{exp-moment-inq-lem-1}, we have the
following exponential moment estimates. Its proof is given in Appendix~\ref{sec5}.

\begin{lem}\label{exp-moment-estimates-lem-1} (1) For any $j=0,1,2$
and $r\in\mathbb R$,
%
\begin{eqnarray}
\label{exp-moment-estimates-lem-1-eq-1}\qquad
&&\limsup_{n\to\infty}
\frac{ b(n)}{ a^2 (b(n))}\nonumber\\[-6pt]\\[-10pt]
&&\qquad{}\times\log E_n \Biggl(\exp \Biggl(\frac{ra^2(b(n))}{ b^2(n)}\sum
_{k\in M_{n\varrho}^{c}} \sum_{l=0}^j
\bigl(\delta_{kl}(n)-E_n\bigl(\delta_{kl}(n)
\bigr)\bigr) \Biggr) \Biggr) \leq0.\nonumber
\end{eqnarray}

(2) If $c_1\in(0,1)$, then for any $j=0,1,2$ and $r\in\mathbb R$,
%
\begin{eqnarray}
\label{exp-moment-estimates-lem-1-eq-2}\qquad
&&
\limsup_{L\to\infty} \limsup
_{n\to\infty} \frac{ b(n)}{ a^2
(b(n))}\nonumber\\[-6pt]\\[-10pt]
&&\qquad{}\times\log E_n \Biggl(\exp
\Biggl(\frac{ra(b(n))}{ b(n)}\sum_{k\in
M_n^{Lc}} \sum
_{l=0}^j\bigl(\delta_{kl}(n)-E_n
\bigl(\delta_{kl}(n)\bigr)\bigr) \Biggr) \Biggr)\leq0\nonumber
\end{eqnarray}
and
%
\begin{eqnarray}
\label{exp-moment-estimates-lem-1-eq-3}\qquad
&&
\limsup_{L\to\infty} \limsup
_{n\to\infty} \frac{ b(n)}{ a^2
(b(n))}\nonumber\\[-6pt]\\[-10pt]
&&\qquad{}\times\log E_n \biggl(\exp
\biggl(\frac{r na(b(n))}{ b(n)}\sum_{k\in
M_n^{Lc}} p_{kn}
\bigl(\delta_{k0}(n)-E_n\bigl(\delta_{k0}(n)
\bigr)\bigr) \biggr) \biggr) \leq0.\nonumber
\end{eqnarray}
\end{lem}

\subsection{\texorpdfstring{The proof of Theorem \protect\ref{main-thm-mdp}: Case \textup{I}}
{The proof of Theorem 2.1: Case I}}\label{sec3.2}

In this subsection, we use the G\"artner--Ellis theorem to show Theorem
\ref{main-thm-mdp} under $c_1\in(0,1)$. The Laplace asymptotic
integral plays a very important role.

By Lemma~\ref{truncation-lem-L}, if $c\in(0,1)$, when $L$ is large enough,
%
\begin{equation}
\label{truncation-bn-eq-1} b^L(n):=E_n
\bigl(F_1^L(n)\bigr) \bigl(1-E_n
\bigl(F_1^L(n)\bigr)/n\bigr)+2E_n
\bigl(F_2^L(n)\bigr)\asymp n
\end{equation}
and
%
\begin{equation}
\label{truncation-bn-eq-2} \lim_{L\to\infty}\limsup
_{n\to\infty}\biggl\llvert \frac
{b^L(n)}{b(n)}-1\biggr\rrvert =0,
\end{equation}
where $ F_j^L(n)=\sum_{k\in M_n^L} \delta_{kj}(n), j\geq1$.
In this subsection, We assume that $L$ is large enough such that
$b^L(n)\asymp n$ and $a(b^L(n))\asymp a(n)$.

The following Laplace asymptotic integral is a key lemma. It will be
proved in Appendix~\ref{sec5}.

\begin{lem}\label{laplace-int-lem}Suppose that conditions (\ref
{moment-condition-1}) and (\ref{moment-condition-2}) hold. If $c_1\in
(0,1)$, then for any $\alpha\in\mathbb R$,
%
\begin{eqnarray}
\label{Laplace-int-lem-eq-0}\qquad
&&
\lim_{n\to\infty} \frac{ b^L(n)}{ a^2 (b^L(n))}\nonumber\\[-6pt]\\[-10pt]
&&\qquad{}\times\log
E_n \biggl(\exp \biggl\{\frac{\alpha a(b^L(n))}{
b^L(n)} \sum
_{k\in M_n^L}\bigl(\delta_{k1}(n)-np_{kn}\delta
_{k0}(n)\bigr) \biggr\} \biggr)= \frac{\alpha^2}{2}.\nonumber
\end{eqnarray}
\end{lem}

\begin{pf*}{Proof of Theorem~\ref{main-thm-mdp} under $c_1\in(0,1)$}
By the G\"artner--Ellis theorem [cf. Theorem 2.3.6 in
\citet{DemZei98}] and Lemma~\ref{laplace-int-lem}, $\{\frac
{1}{a(b^L(n))}\sum_{k\in M_n^L}(\delta_{k1}(n)-np_{kn}\delta
_{k0}(n)),n\geq1\}$ satisfies a large deviation principle with speed
$\frac{a^2(b^L(n))}{b^L( n)}$ and with rate function $I(x)=\frac
{x^2}{2}$. By Lemma~\ref{exp-eqv-lemma-1}, we only need to check
%
\begin{eqnarray}
\label{main-thm-case-I-eq-0}\qquad
&&
\limsup_{L\to\infty}\limsup
_{n\to\infty} \frac{ b(n)}{ a^2
(b(n))}
\nonumber\\
&&\quad{}\times\log P_n \biggl( \biggl\llvert \frac{1}{a(b^L(n))}\sum
_{k\in
M_n^L}\bigl(\delta _{k1}(n)-np_{kn}
\delta_{k0}(n)\bigr)-\frac{n(\hat
{Q}_{n}-Q_n)}{a(b(n))}\biggr\rrvert \geq\epsilon
\biggr)\\
&&\qquad=-\infty.
\nonumber
\end{eqnarray}

It is obvious that
%
\begin{eqnarray}
\label{main-thm-case-I-eq-1}
&& P_n \biggl( \biggl\llvert
\frac{1}{a(b^L(n))}\sum_{k\in M_n^L}\bigl(\delta
_{k1}(n)-np_{kn}\delta_{k0}(n)\bigr)-
\frac{n(\hat
{Q}_{n}-Q_n)}{a(b(n))}\biggr\rrvert \geq\epsilon \biggr)
\nonumber
\\
&&\qquad\leq P_n \biggl( \biggl\llvert \frac{a(b^L(n))-a(b(n))}{a(b^L(n))a(b(n))}\sum
_{k\in M_n^L}\bigl(\delta_{k1}(n)-np_{kn}
\delta_{k0}(n)\bigr)\biggr\rrvert \geq \epsilon /2 \biggr)
\\
&&\qquad\quad{}+P_n \biggl( \frac{1}{a(b(n))}\biggl\llvert \sum
_{k\in M_n^{Lc}}\bigl(np_{kn} \delta _{k0}(n)-
\delta_{k1}(n)\bigr)\biggr\rrvert \geq\epsilon/2 \biggr).
\nonumber
\end{eqnarray}
From (\ref{truncation-bn-eq-2}) and $\{\frac{1}{a(b^L(n))}\sum_{k\in
M_n^L}(\delta_{k1}(n)-np_{kn}\delta_{k0}(n)),n\geq1\}$ satisfies the
large deviation principle, we obtain that for any $\epsilon>0$,
%
\begin{eqnarray}
\label{main-thm-case-I-eq-2}\qquad
&&\limsup_{L\to\infty}\limsup
_{n\to\infty} \frac{ b^L(n)}{ a^2
(b^L(n))}
\nonumber\\
&&\quad{}\times\log P_n \biggl( \biggl\llvert \frac
{a(b^L(n))-a(b(n))}{a(b^L(n))a(b(n))}\sum
_{k\in M_n^L}\bigl(\delta _{k1}(n)-np_{kn}
\delta_{k0}(n)\bigr)\biggr\rrvert \geq\epsilon \biggr)\\
&&\qquad=-\infty.
\nonumber
\end{eqnarray}
By Lemma~\ref{exp-moment-estimates-lem-1} and Chebyshev's inequality,
we have that for any $\epsilon>0$,
\begin{eqnarray*}
&&
\limsup_{L\to\infty}\limsup_{n\to\infty}
\frac{ b(n)}{ a^2
(b(n))}\log P_n \biggl( \frac{n}{a(b(n))}\biggl\llvert
\sum_{k\in M_n^{Lc}}p_{kn} \bigl(\delta
_{k0}(n)-E_n\bigl(\delta_{k0}(n)\bigr)\bigr)
\biggr\rrvert \geq\epsilon \biggr)\\
&&\qquad=-\infty
\end{eqnarray*}
and for $j=0,1$,
\begin{eqnarray*}
&&
\limsup_{L\to\infty}\limsup_{n\to\infty}
\frac{ b(n)}{ a^2
(b(n))}\log P_n \Biggl( \frac{1}{a(b(n))}\Biggl\llvert
\sum_{k\in M_n^{Lc}} \sum_{l=0}^j
\bigl(\delta_{kl}(n) - E_n\bigl(\delta_{kl}(n)
\bigr)\bigr)\Biggr\rrvert \geq\epsilon \Biggr)\\
&&\qquad=-\infty,
\end{eqnarray*}
which implies that for any $\epsilon>0$,
%
\begin{eqnarray}
\label{main-thm-case-I-eq-3}
&&
\limsup_{L\to\infty}\limsup
_{n\to\infty} \frac{ b(n)}{ a^2
(b(n))}\nonumber\\
&&\quad{}\times\log P_n \biggl(
\frac{1}{a(b(n))} \biggl|\sum_{k\in M_n^{Lc}}\bigl(np_{kn}
\delta _{k0}(n)- \delta_{k1}(n)\bigr) \biggr|\geq\epsilon
\biggr)\\
&&\qquad=-\infty.\nonumber
\end{eqnarray}
Now, (\ref{main-thm-case-I-eq-0}) follows from (\ref
{main-thm-case-I-eq-2}) and (\ref{main-thm-case-I-eq-3}). Therefore,
the conclusion of Theorem~\ref{main-thm-mdp} holds under $c_1\in(0,1)$.
\end{pf*}

\subsection{\texorpdfstring{The proof of Theorem \protect\ref{main-thm-mdp}: Case \textup{II}}
{The proof of Theorem 2.1: Case II}}\label{sec3.3}

In this subsection, we show Theorem~\ref{main-thm-mdp} under $c_1=0$.
In this case, since $\{np_{in},i\geq1,n\geq1\}$ cannot be truncated
as a uniformly bounded sequence, the asymptotic analysis techniques in
the first case cannot be used.
The proof of this case is based
on the Poisson approximation [cf. \citet{ZhaZha09}] and the
truncation method.

Let first us introduce the Poissonization defined by \citet{ZhaZha09}. Define
%
\begin{equation}
\label{def-xi-n} \xi_n=\sum_{k=1}^\infty
\bigl(\delta_{k1}(n)-np_{kn}\delta_{k0}(n)
\bigr)=n(\hat{Q}_n-Q_n).
\end{equation}
Let $N_\lambda$ be a Poisson process independent of $\{{\mathbf
X}(m),m\geq1\}$ with $E_n(N_\lambda)=\lambda$. Define the
Poissonization $\zeta_{\lambda n} $ of $\xi_n$ as follows:
%
\begin{equation}
\label{def-zeta-n} \zeta_{\lambda n} = \sum_{k=1}^\infty
Y_{k\lambda n}\qquad\mbox{where } Y_{k\lambda
n}=\delta_{k1}(N_\lambda)-
\lambda p_{kn}\delta_{k0}(N_\lambda).
\end{equation}

Under probability $P_n$, $ X_k(N_\lambda), k\geq1 $ are independent
Poisson variables with means $\lambda p_{kn}$, so that $ Y_{k\lambda
n},k\geq1 $ are independent zero-mean variables with variance
$\sigma_{k \lambda n}^2:=\lambda p_{kn}e^{-\lambda p_{kn}}+(\lambda
p_{kn})^2 e^{-\lambda p_{kn}}$.
Then the Poissonization $\{\zeta_{n n},n\geq1\} $ satisfies the
following moderate deviation principle.

\begin{lem}\label{mdp-poisson-approx}
Let conditions (\ref{moment-condition-1}), (\ref{moment-condition-2})
and (\ref{mdp-Lindeberg-condition}) hold. Then
$\{\frac{\zeta_{nn}}{a(s_n^2)},\break n\geq
1\}$ satisfies a large deviation principle with speed
$\frac{a^2(s_n^2)}{ s_n^2 }$ and with rate function \mbox{$I(x)=\frac{x^2}{2}$}.
\end{lem}

\begin{pf}
For any $\alpha\in\mathbb R$,
\begin{eqnarray*}
&&
E_n \biggl(\exp \biggl\{\frac{\alpha a (s_n^2)}{
s_n^2}\zeta_{nn}
\biggr\} \biggr) \\
&&\qquad= \prod_{k=1}^\infty
E_n \biggl(\exp \biggl\{\frac{\alpha a
(s_n^2)}{ s_n^2} (I_{\{X_k(N_{\lambda_n})=1\}}
-np_{kn} I_{\{X_k(N_{\lambda_n})=0\}} ) \biggr\} \biggr)
\\
&&\qquad= \prod_{k=1}^\infty \biggl(
\bigl(1-e^{-np_{kn}}-np_{kn}e^{-np_{kn}} \bigr)+np_{kn}e^{-np_{kn}}
\exp \biggl\{\frac{\alpha a (s_n^2)}{
s_n^2} \biggr\}
\\
&&\qquad\quad\hspace*{132pt}{} + e^{-np_{kn}}\exp \biggl\{\frac{-\alpha a
(s_n^2)}{ s_n^2} np_{kn} \biggr\}
\biggr).
\end{eqnarray*}

For any $\epsilon\in(0,1/2 ]$ such that $|\alpha|\epsilon<1/2$, for
$n$ large enough, we can write
\begin{eqnarray*}
&& 1-e^{-np_{kn}}-np_{kn}e^{-np_{kn}} \\
&&\quad{}+np_{kn}e^{-np_{kn}}
\exp \biggl\{\frac{\alpha a
(s_n^2)}{ s_n^2} \biggr\} + e^{-np_{kn}}\exp \biggl\{
\frac{-\alpha a (s_n^2)}{
s_n^2} np_{kn} \biggr\}
\\
&&\qquad= 1 +\frac{1}{2} \biggl(\frac{\alpha a (s_n^2) }{ s_n^2} \biggr)^2
\bigl(np_{kn}+(np_{kn})^2 \bigr)e^{-np_{kn}}
+o \biggl( \biggl(\frac{ a (s_n^2) }{ s_n^2} \biggr)^2 \biggr)np_{kn}e^{-np_{kn}}
\\
&&\qquad\quad{}+ \biggl(\frac{\alpha a (s_n^2) }{ s_n^2}np_{kn}-\frac{1}{2} \biggl(
\frac
{\alpha a (s_n^2) }{ s_n^2} \biggr)^2(np_{kn})^2
\biggr)e^{-np_{kn}}I_{\{np_{kn}>\epsilon
s_n^2/a(s_n^2)\}}
\\
&&\qquad\quad{}+O \biggl( \biggl(\frac{ a (s_n^2) }{ s_n^2} \biggr)^3 \biggr)
(np_{kn})^3e^{-np_{kn}}I_{\{np_{kn}\leq\epsilon
s_n^2/a(s_n^2)\}}
\\
&&\qquad\quad{}+e^{-np_{kn}} \biggl(\exp \biggl\{\frac{-\alpha a (s_n^2)}{
s_n^2} np_{kn}
\biggr\}-1 \biggr)I_{\{np_{kn}>\epsilon s_n^2/a(s_n^2)\}}.
\end{eqnarray*}
By (\ref{mdp-Lindeberg-condition}),
\begin{eqnarray*}
&&
\frac{1}{a(s_n^2)} \sum_{k=1}^\infty
np_{kn} e^{-np_{kn}}I_{\{
np_{kn}>\epsilon s_n^2/a(s_n^2)\}}\\
&&\qquad\leq\frac{1}{\epsilon s_n^2} \sum
_{k=1}^\infty(np_{kn})^2
e^{-np_{kn}}I_{\{np_{kn}>\epsilon
s_n^2/a(s_n^2)\}}\to0
\end{eqnarray*}
and
\begin{eqnarray*}
&&\frac{s_n^2}{a^2(s_n^2)} \sum_{k=1}^\infty
e^{-np_{kn}} \biggl(\exp \biggl\{\frac{-\alpha a (s_n^2)}{
s_n^2} np_{kn} \biggr
\}-1 \biggr)I_{\{np_{kn}>\epsilon s_n^2/a(s_n^2)\}
}
\\
&&\qquad\leq\frac{2 }{\epsilon^2 s_n^2} \sum_{k=1}^\infty(np_{kn})^2
e^{-np_{kn}(1-|\alpha| a (s_n^2)/s_n^2)}I_{\{np_{kn}>\epsilon
s_n^2/a(s_n^2)\}} \to0.
\end{eqnarray*}
Therefore, by $
\frac{ a (s_n^2) }{ s_n^2}\frac{ 1}{ s_n^2} \sum_{k=1}^\infty
(np_{kn})^3e^{-np_{kn}}I_{\{np_{kn}\leq\epsilon s_n^2/a(s_n^2)\}}\leq
\varepsilon\to0$ as $\epsilon\to0$, we have that
\[
\lim_{n\to\infty} \frac{ s_n^2}{ a^2 (s_n^2)}\log E_n \biggl(
\exp \biggl\{ \frac{\alpha a (s_n^2)}{
s_n^2}\zeta_{nn} \biggr\} \biggr) =
\frac{\alpha^2}{2},
\]
which implies the conclusion of the lemma by the G\"artner--Ellis
theorem; cf. Theorem 2.3.6 in \citet{DemZei98}.
\end{pf}

By Lemmas~\ref{mdp-poisson-approx} and~\ref{exp-eqv-lemma-basic},
we need the
following exponential approximation:
for any
$\varepsilon>0$,
%
\begin{equation}
\label{exp-eqv-lemma-2-eq-1} \lim_{n\to\infty}\frac{b( n)}{a^2(b(n))}\log
P_n\bigl(|\xi_n-\zeta_{nn}|>\varepsilon a
\bigl(b(n)\bigr)\bigr)=-\infty.
\end{equation}

Let us first give a maximal exponential estimate. Its proof is
postponed to Appendix~\ref{sec5}.

\begin{lem}\label{exp-eqv-lemma-1}
Let conditions (\ref{moment-condition-1}), (\ref{moment-condition-2})
and (\ref{mdp-Lindeberg-condition}) hold, and let
$c_1=0$. For
any $M\geq1$ fixed, set $\lambda_n=n- M a(b(n))\sqrt{\frac{n}{b(n)}}$,
$\Delta_n=2M
a(b(n))\sqrt{\frac{n}{b(n)}}$. Then for any $\varepsilon>0$,
%
\begin{equation}
\label{exp-eqv-lemma-1-eq-1}\qquad \lim_{n\to\infty}\frac{b( n)}{a^2(b(n))}\log
P_n \Bigl(\max_{t\in
[\lambda_n, \lambda_n+ \Delta_n ]}|\zeta_{\lambda_n n}-
\zeta_{t
n}|\geq\varepsilon a\bigl(b(n)\bigr) \Bigr)=-\infty.
\end{equation}
\end{lem}

\begin{pf*}{Proof of Theorem~\ref{main-thm-mdp} under $c_1=0$}
By Lemmas~\ref{mdp-poisson-approx} and~\ref{exp-eqv-lemma-basic},
we only need to prove (\ref{exp-eqv-lemma-2-eq-1}).
Set $t_n=\inf\{\lambda;N_\lambda=n\}$. Then $t_n$ has gamma$(n,1)$
distribution and $\xi_n-\zeta_{t_nn}=(t_n-n)\sum_{k=1}^\infty
p_{kn}\delta_{k0}(n)$. Therefore,\vadjust{\goodbreak} for any $\varepsilon>0$ and any
$M\geq1$,
\begin{eqnarray*}
&& P_n \bigl(|\xi_n-\zeta_{nn}|\geq
\varepsilon a\bigl(b(n)\bigr) \bigr)
\\
&&\qquad\leq P_n \biggl(|t_n-n|\geq M a\bigl(b(n)\bigr)\sqrt{
\frac{n}{b(n)}} \biggr)+P_n \Biggl( \sum
_{k=1}^\infty p_{kn}\delta_{k0}(n)
\geq\frac
{\varepsilon}{2M}\sqrt{\frac{b(n)}{n}} \Biggr)
\\
&&\qquad\quad{}+ P_n \biggl(\max_{t\in[n- \Delta_n/2, n+ \Delta_n/2 ]}|\zeta _n-
\zeta _{t n}|\geq\frac{\varepsilon a(b(n))}{2} \biggr).
\end{eqnarray*}
By Lemma~\ref{exp-eqv-lemma-1},
%
\begin{equation}
\label{exp-eqv-lemma-2-eq-2} \lim_{n\to\infty}\frac{b( n)}{a^2(b(n))}\log
P_n \biggl(\max_{t\in[n-
\Delta_n/2, n+ \Delta_n/2 ]}|\zeta_n-
\zeta_{t n}|\geq \frac{\varepsilon a(b(n))}{2}
\biggr)=-\infty.\hspace*{-34pt}
\end{equation}
By Chebyshev's inequality, it is easy to get that
%
\begin{equation}
\label{exp-eqv-lemma-2-eq-3}\qquad \lim_{M\to\infty}\limsup
_{n\to\infty}\frac{b( n)}{a^2(b(n))}\log P_n
\biggl(|t_n-n|\geq M a\bigl(b(n)\bigr)\sqrt{\frac{n}{b(n)}}
\biggr)=-\infty.
\end{equation}
Therefore, we only need to prove that
%
\begin{equation}
\label{exp-eqv-lemma-2-eq-4} \limsup_{n\to\infty}\frac{b( n)}{a^2(b(n))}
\log P_n \Biggl( \sum_{k=1}^\infty
p_{kn}\delta_{k0}(n)\geq \frac{\varepsilon}{2M} \sqrt{
\frac{b(n)}{n}} \Biggr)=-\infty.
\end{equation}
It is sufficient that for any $r>0$,
%
\begin{equation}
\label{exp-eqv-lemma-2-eq-5}\qquad \lim_{n\to\infty}\frac{ b(n)}{ a^2(b(n))}\log
E_n \Biggl(\exp \Biggl(\frac{r a^2(b(n))}{ b(n)}\sqrt{\frac
{n}{b(n)}}
\sum_{k=1}^\infty p_{kn}
\delta_{k0}(n) \Biggr) \Biggr)=0.
\end{equation}
In fact, by Lemma~\ref{exp-moment-inq-lem-1}, we can get that for any $r>0$,
\begin{eqnarray*}
&&\log E_n \Biggl(\exp \Biggl(\frac{r a^2(b(n))}{ b(n)}\sqrt{
\frac
{n}{b(n)}}\sum_{k=1}^\infty
p_{kn}\delta_{k0}(n) \Biggr) \Biggr)
\\
&&\qquad\leq \frac{2 r a^2(b(n))}{b(n)}\sqrt{\frac{n}{b(n)}} \\
&&\qquad\quad{}\times\sum
_{k=1}^\infty \biggl( p_{kn} e^{-n p_{kn}}
+ p_{kn}\exp \biggl\{-n \biggl(1-\frac{r a^2(b(n))}{ n b(n)}\sqrt{
\frac
{n}{b(n)}} \biggr) p_{kn} \biggr\} \biggr)
\\
&&\qquad\leq\frac{ r a^2(b(n)) }{ b(n)}\sqrt{\frac{n}{b(n)}} \biggl( \frac{2
s_n^2}{n}+
\frac{
s_{\lambda_n n}^2}{\lambda_n} \biggr),
\end{eqnarray*}
where $\lambda_n=n (1-\frac{r a^2(b(n))}{ n b(n)} \sqrt{\frac
{n}{b(n)}} )$,
which implies that (\ref{exp-eqv-lemma-2-eq-1}) holds.
\end{pf*}

\subsection{\texorpdfstring{Proof of Theorem \protect\ref{main-thm-mdp-self}}{Proof of Theorem 2.2}}\label{sec3.4}

By the comparison method in large deviations [cf. Theorem 4.2.13 in
\citet{DemZei98}], in order to obtain Theorem \ref
{main-thm-mdp-self}, we need the following lemma.

\begin{lem}\label{main-thm-self-lem} For any $\varepsilon>0$, for $j=1,2$,
%
\begin{equation}
\label{main-thm-self-lem-eq-1} \limsup_{n\to\infty}\frac{b( n)}{a^2(b(n))}
\log P_n \bigl( \bigl|F_j(n)-E_n
\bigl(F_j(n)\bigr)\bigr|\geq\varepsilon b(n) \bigr)=-\infty.
\end{equation}
\end{lem}

\begin{pf}
By (\ref{exp-moment-estimates-lem-1-eq-1}),
for $j=1,2$, for any $\varrho>0$ and $\varepsilon>0$,
%
\begin{equation}
\label{main-thm-self-lem-eq-2}\quad \limsup_{n\to\infty}
\frac{ b(n)}{ a^2 (b(n))} \log P_n \biggl( \frac
{1}{b(n)} \biggl|\sum
_{k\in M_{n\varrho}^{c}}\bigl( \delta_{kj}(n)-
E_n\bigl(\delta _{kj}(n)\bigr)\bigr) \biggr|\geq\varepsilon
\biggr)=-\infty.\hspace*{-36pt}
\end{equation}
Therefore, by Lemma~\ref{truncation-lem-rho}, it suffices to show that
%
\begin{eqnarray}
\label{main-thm-self-lem-eq-3}
&&
\limsup_{\varrho\to0}\limsup
_{n\to\infty} \frac{ b(n)}{ a^2
(b(n))}
\nonumber\\[-2pt]
&&\quad{}\times\log P_n \biggl( \frac{1}{b(n)} \biggl|\sum
_{k\in M_{n\varrho
}} \biggl(\delta_{kj}(n)- \frac{1}{j!}(np_{kn})^j
e^{-np_{kn}} \biggr) \biggr|\geq \varepsilon \biggr)\\[-2pt]
&&\qquad=-\infty.
\nonumber
\end{eqnarray}

Now, let us show (\ref{main-thm-self-lem-eq-3}). Using the partial
inversion formula for characteristic function due to \citet{Bar}
[see also \citet{Hol79}, \citet{Est83}], for any $r\in\mathbb R$,
\begin{eqnarray*}
\hspace*{-4pt}&&E_n \Biggl(\exp \Biggl\{ r\sum_{k=1}^\infty
\biggl(\delta _{kj}(n)-\frac
{1}{j!}(np_{kn})^j
e^{-np_{kn}} \biggr) \Biggr\} \Biggr)
\\[-2pt]
\hspace*{-4pt}&&\qquad= \frac{n!}{2\pi n^ne^{-n}}\int_{-\pi}^\pi\prod
_{k\in M_{n\varrho}^c} E_n \bigl(\exp \bigl\{iu
\bigl(Y_k(n)-np_{kn}\bigr) \bigr\} \bigr)
\\[-2pt]
\hspace*{-4pt}&&\hspace*{60pt}\qquad\quad{}\times \prod_{k\in M_{n\varrho}} E_n \biggl(\exp
\biggl\{iu \bigl(Y_k(n)-np_{kn}\bigr)\\[-2pt]
\hspace*{-4pt}&&\qquad\hspace*{147pt}{} +r
\biggl(I_{\{Y_k(n) = j\}}-\frac{1}{j!}(np_{kn})^j
e^{-np_{kn}} \biggr) \biggr\} \biggr)\,du,
\end{eqnarray*}
where $Y_k(n),k\geq1$ are independent random variables and $Y_k(n)$ is
Poisson distributed with mean $np_{kn}$.
Let $\gamma_k(u)$ be defined as in the proof of Lemma~\ref
{laplace-int-lem}, that is,
$
\gamma_k(u)= \exp \{np_{kn}(e^{iu }-1-iu )  \}.
$
Set
\begin{eqnarray*}
\hspace*{-3pt}&&\vartheta_k(u,\alpha)
\\
\hspace*{-3pt}&&\quad= \biggl( \exp \bigl\{iju- np_{kn}\bigl(e^{iu }-1\bigr)
\bigr\} \biggl(\exp \biggl\{ \frac{\alpha a^2(b(n))}{b^2(n)} \biggr\}-1 \biggr)
\frac
{1}{j!}(np_{kn})^j e^{-np_{kn}} +1 \biggr)
\\
\hspace*{-3pt}&&\qquad\times\exp \biggl\{-\frac{\alpha a^2(b(n))}{b^2(n)}\frac
{1}{j!}(np_{kn})^j
e^{-np_{kn}} \biggr\}.
\end{eqnarray*}
Then for any $\alpha\in\mathbb R$,
\begin{eqnarray*}
&&
E_n \biggl(\exp \biggl\{ \frac{\alpha a^2(b(n))}{b^2(n)} \sum
_{k\in
M_{n\varrho}} \biggl(\delta_{kj}(n)-\frac{1}{j!}(np_{kn})^j
e^{-np_{kn}} \biggr) \biggr\} \biggr)
\\
&&\qquad= \frac{n!}{2\pi n^ne^{-n}}\int_{-\pi}^\pi
e^{n(e^{i u}-1-iu)} \prod_{k\in M_{n\varrho}}\vartheta_k(u,
\alpha)\,du.
\end{eqnarray*}
Set\vspace*{1pt} $\tau(n)=\sqrt{a(b(n)/b(n)}$. Then
$\frac{\log n}{n\tau^2(n)}=\frac{b(n)\log n}{n a(b(n))}\leq\frac
{\sqrt {b(n)}\log n}{n}\to0$, and noting that
$
\sum_{k\in M_{n\varrho}}np_{kn}\leq n$, $\sum_{k\in M_{n\varrho
}}(np_{kn})^2\leq\varrho n b(n)/a(b(n)),
$
we obtain\vspace*{1pt} that for $\varrho$ small enough,
\begin{eqnarray*}
&&
\frac{ b(n)}{ a^2 (b(n))}\log \biggl( n ^{1/2}\sup_{ |u| \in[\tau
(n),\pi]}
\biggl\llvert e^{n(e^{i u}-1-iu)} \prod_{k\in M_{n\varrho}} \vartheta
_k(u,\alpha)\biggr\rrvert \biggr)
\\
&&\qquad\leq - \frac{ b(n)n \tau^2(n)}{a^2(b(n))} \biggl(1+O \biggl(\frac
{\log
n}{ n\tau^2(n)} \biggr)+O(
\varrho) \biggr) \to-\infty.
\end{eqnarray*}
Since $\sup_{u\in[-\tau(n),\tau(n)]}\sup_{k\in M_{n\varrho
}}|np_{kn}(1-\cos u)|\leq\varrho$, on $[-\tau(n),\tau(n)]$,
\begin{eqnarray*}
\hspace*{-3pt}&&\biggl\llvert e^{n(e^{i u}-1-iu)} \prod_{k\in M_{n\varrho}}
\vartheta _k(u,\alpha) \biggr\rrvert
\\
\hspace*{-3pt}&&\quad=\exp \biggl\{ \biggl(O(\varrho) +O \biggl(\frac{\alpha^2
a^2(b(n))}{b^2(n)} \biggr) \biggr)
\frac{a^2(b(n))}{b(n)} \biggr\}\exp \biggl\{-\frac{n}{2}u^2
\bigl(1+ O(u)o(1) \bigr) \biggr\}.
\end{eqnarray*}
Thus
\[
\limsup_{\varrho\to0}\limsup_{n\to\infty}
\frac{ b(n)}{ a^2
(b(n))}\log\biggl\llvert \int_{-\tau(n)}^{\tau(n)}
n^{1/2} e^{n(e^{i
u}-1-iu)} \prod_{k\in M_{n\varrho}}
\vartheta_k(u,\alpha) \,du\biggr\rrvert =0
\]
and so
\[
\limsup_{\varrho\to0} \limsup_{n\to\infty}
\frac{ b(n)}{ a^2
(b(n))}\log E_n \biggl(\!\exp \biggl\{ r\!\sum
_{k\in M_{n\varrho}}\! \biggl(\!\delta _{kj}(n)-\frac{1}{j!}(np_{kn})^j
e^{-np_{kn}} \!\biggr)\! \biggr\}\! \biggr)\leq0.
\]
This yields that (\ref{main-thm-self-lem-eq-3}) holds.
\end{pf}

\begin{pf*}{Proof of Theorem~\ref{main-thm-mdp-self}}
By Lemma~\ref{main-thm-self-lem}, for any $\varepsilon>0$,
\[
\limsup_{n\to\infty}\frac{b( n)}{a^2(b(n))}\log P_n \biggl(
\biggl\llvert \frac
{b(n)}{F_1(n)(1-F_1(n)/n)+2F_2(n)}-1\biggr\rrvert \geq\varepsilon \biggr)=-\infty.
\]
Now, by
\begin{eqnarray*}
&&\biggl\llvert \frac{\sqrt{b(n)}n(\hat{Q}_{n}-Q_n)}{a(b(n)) \sqrt {F_1(n)(1-F_1(n)/n)+2F_2(n) }}- \frac{ n(\hat
{Q}_{n}-Q_n)}{a(b(n))}\biggr\rrvert
\\
&&\qquad= \biggl\llvert \frac{ n(\hat{Q}_{n}-Q_n)}{a(b(n))}\biggr\rrvert \biggl\llvert \sqrt {
\frac
{b(n)} {F_1(n)(1-F_1(n)/n)+2F_2(n)}}-1\biggr\rrvert,
\end{eqnarray*}
and the elementary inequality $|x-1|=|x^2-1|/|x+1|\leq|x^2-1|$ for all
$x\geq0$,
we obtain that
%
\begin{eqnarray}
\label{main-thm-self-proof-eq-2}\qquad
&&\limsup_{n\to\infty}
\frac{b( n)}{a^2(b(n))}\log P_n \biggl( \biggl|\frac
{\sqrt{b(n)}n(\hat{Q}_{n}-Q_n)}{a(b(n)) \sqrt {F_1(n)(1-F_1(n)/n)+2F_2(n) }}
\nonumber\\
&&\qquad\hspace*{200.8pt}{} - \frac{ n(\hat{Q}_{n}-Q_n)}{a(b(n))} \biggr|\geq \varepsilon \biggr)\\
&&\qquad=-\infty.
\nonumber
\end{eqnarray}
Therefore, the conclusion of the theorem follows from Lemma \ref
{exp-eqv-lemma-basic} or Theorem~4.2.13 in \citet{DemZei98}.
\end{pf*}

\begin{appendix}\label{app}
\section{Some concepts of large deviations}\label{sec4}

For the sake convenience, let us introduce some notions in large
deviations [\citet{DemZei98}].
Let $({\mathcal X},\rho)$ be a metric space.
Let $(\Omega_n, \mathcal F_n, P_n)$, $n\geq1$ be a sequence of
probability spaces and let $\{\eta_n, n\geq1\}$ be a sequence of
measurable maps
from $\Omega_n$ to ${\mathcal X}$. Let $\{\lambda_n,n\geq1\}$ be a sequence
of positive numbers tending to $+\infty$, and let $I\dvtx {\mathcal X} \to[0,
+\infty]$ be inf-compact; that is, $[I\le L]$ is compact for any
$L\in\mathbb R$. Then $\{\eta_n,n\geq1\}$ is said to satisfy a large
deviation
principle (LDP) with speed $\lambda_n$ and with rate function $I$, if
for any open measurable subset $G$ of ${\mathcal X}$,
%
\begin{equation}
\label{LLD-def} \liminf_{n\rightarrow\infty}\frac1{\lambda_n}
\log{P_n}(\eta _n\in G) \ge- \inf_{x \in G}I(x)
\end{equation}
and for any closed measurable subset $F$ of ${\mathcal X}$,
%
\begin{equation}
\label{ULD-def} \limsup_{n\rightarrow\infty} \frac1{\lambda_n}
\log {P_{n}}(\eta_n\in F) \le-\inf_{x\in F}
I(x).
\end{equation}

\begin{rmk}
Assume that $\{\eta_n,n\geq1\}$ satisfies $\eta_n\to\mu$ in law
and a
fluctuation theorem such as central limit theorem, that is, there
exists a sequence $l_n\to\infty$ such
that $l_n(\eta_n-\mu)\rightarrow\eta$ in law, where $\mu$ is a
constant and $\eta$ is a nontrivial random variable. Usually, $\{\eta
_n, n\geq1\}$ is said to satisfy a moderate
deviation principle (MDP) if $\{r_n(\eta_n-\mu),n \geq1\}$ satisfies
a large
deviation principle, where $r_n$ is an intermediate scale between $1$
and $l_n$, that is, $r_n\to\infty$ and $r_n/l_n\to0$.
\end{rmk}

In this paper, the following exponential approximation lemma is
required. It is slightly different from Theorem 4.2.16 in \citet{DemZei98}.

\begin{lem}\label{exp-eqv-lemma-basic}
Let $\{\eta_n, n\geq1\}$ and $\{\eta_n^L, n\geq1\}$, $L\geq1$ be
sequences of measurable maps
from $\Omega_n$ to ${\mathcal X}$. Assume that for each $L\geq1$, $\{
\eta_n^L,n\geq1\}$ satisfies a LDP with speed $\lambda^L_n$ and with
rate function $I$. If
%
\begin{equation}
\label{exp-eqv-lemma-basic-eq-1} \lim_{L\to\infty}\limsup
_{n\rightarrow\infty}\biggl\llvert \frac
{\lambda
_n^L}{\lambda_n}-1\biggr\rrvert =0
\end{equation}
and for any $\epsilon>0$,
%
\begin{equation}
\label{exp-eqv-lemma-basic-eq-2} \lim_{L\to\infty}\limsup
_{n\rightarrow\infty}\frac1{\lambda_n} \log {P_n}
\bigl(\rho\bigl(\eta_n,\eta_n^L\bigr)\geq
\epsilon\bigr)=-\infty,
\end{equation}
the $\{\eta_n, n\geq1\}$ satisfies a LDP with speed $\lambda_n$ and
with rate function $I$.
\end{lem}
\begin{pf} Set $I(A)=\inf_{x\in A} I(x)$.
For any closed subset $F$,
\[
P (\eta_n\in F )\leq P \bigl(\eta_n^L\in
F^\epsilon \bigr)+P \bigl(\rho\bigl(\eta_n,
\eta_n^L\bigr)\geq\epsilon \bigr),
\]
where $F^\epsilon=\{y\in\mathcal X;\inf_{x\in F}\rho(y,x)<\epsilon
\}$.
By (\ref{exp-eqv-lemma-basic-eq-2}),
\[
P \bigl(\rho\bigl(\eta_n,\eta
_n^L\bigr)\geq
\epsilon \bigr)\leq e^{-\lambda_n(I(F^\epsilon)+1)}
\]
for large $n$ and
$L$. Therefore, for large $n$ and $L$
\[
P \bigl(\eta_n^L\in F^\epsilon \bigr)+P \bigl(
\rho\bigl(\eta_n,\eta _n^L\bigr)\geq
\epsilon \bigr) \leq e^{-\lambda_n(I(F^\epsilon)+o(1))}+e^{-\lambda_n(I(F^\epsilon)+1)}
\]
and so
\[
\limsup_{n\to\infty}\frac{1}{\lambda_n}\log P (\eta_n
\in F )\leq- I\bigl(F^\epsilon\bigr)\to-\inf_{x\in F}I(x).
\]

The argument for open sets is similar and is omitted.
\end{pf}

\section{\texorpdfstring{Proofs of Lemmas \lowercase{\protect\ref{exp-moment-estimates-lem-1}},
\lowercase{\protect\ref{laplace-int-lem}} and \lowercase{\protect\ref{exp-eqv-lemma-1}}}
{Proofs of Lemmas 3.6, 3.7 and 3.9}}\label{sec5}

In this Appendix, we give the proofs of several technique lemmas. The
proofs of Lemmas~\ref{exp-moment-estimates-lem-1} and \ref
{exp-eqv-lemma-1} are based some exponential moment inequalities for
negatively associated random variables and martingales. The refined
asymptotic analysis techniques play a basic role in the proof of Lemma
\ref{laplace-int-lem}.

\begin{pf*}{Proof of Lemma~\ref{exp-moment-estimates-lem-1}} (1) By
Lemma~\ref{exp-moment-inq-lem-1}, we have that for any $r\in\mathbb R$,
and $j=0,1,2$,
\begin{eqnarray*}
\hspace*{-2pt}&& \log E_n \Biggl(\exp \Biggl(\frac{ra^2(b(n))}{ b^2(n)}\sum
_{k\in
M_{n\varrho}^{c}} \sum_{l=0}^j
\bigl(\delta_{kl}(n)-E_n\bigl(\delta _{kl}(n)
\bigr)\bigr) \Biggr) \Biggr)
\\
\hspace*{-2pt}&&\qquad\leq\sum_{k\in M_{n\varrho}^{c}} \Biggl( \log \Biggl( \biggl(\exp
\biggl\{ \frac{ra^2(b(n))}{ b^2(n)} \biggr\}-1 \biggr)\sum_{l=0}^j
\frac
{n!}{(n-l)!l!} p_{kn}^l (1-p_{kn})^{n-l}+1
\Biggr)
\\
\hspace*{-2pt}&&\hspace*{125pt}\qquad\quad{} -\frac{ra^2(b(n))}{ b^2(n)}\sum_{l=0}^j
\frac{n!}{(n-l)!l!} p_{kn}^l (1-p_{kn})^{n-l}
\Biggr)
\\[-4pt]
\hspace*{-2pt}&&\qquad\leq\frac{a(b(n))}{b(n)}\frac{r^2a^2(b(n))}{b(n)} \\[-4pt]
\hspace*{-2pt}&&\qquad\quad{}\times\frac{1}{b(n)} \sum
_{k\in M_{n\varrho}^{c}} \Biggl( \frac{2}{\varrho} np_{kn}e^{-n p_{kn}}+
\sum_{l=1}^j\frac{n!}{(n-l)!l!}
p_{kn}^l e^{-(n-l)p_{kn}} \Biggr).
\end{eqnarray*}
Therefore, (\ref{exp-moment-estimates-lem-1-eq-1}) holds.

(2) Similarly to the proof of (\ref{exp-moment-estimates-lem-1-eq-1}),
we also have that
\begin{eqnarray*}
&& \log E_n \Biggl(\exp \Biggl(\frac{ra(b(n))}{ b(n)}\sum
_{k\in
M_n^{Lc}}\sum_{l=0}^j
\bigl(\delta_{kl}(n)-E_n\bigl(\delta_{kl}(n)
\bigr)\bigr) \Biggr) \Biggr)
\\[-4pt]
&&\qquad\leq\frac{r^2a^2(b(n))}{b(n)} \frac{1}{b(n)}\sum_{k\in
M_n^{Lc}}
\Biggl( \frac{2}{L } np_{kn}e^{-n p_{kn}}+ \sum
_{l=1}^j\frac{n!}{(n-l)!l!} p_{kn}^l
e^{-(n-l)p_{kn}} \Biggr)\\[-4pt]
&&\qquad\to0.
\end{eqnarray*}
Finally, let us prove (\ref{exp-moment-estimates-lem-1-eq-3}). By Lemma
\ref{exp-moment-inq-lem-1}, for any $r\not=0$,
\begin{eqnarray*}
&& \log E_n \biggl(\exp \biggl(\frac{r na(b(n))}{ b(n)}\sum
_{k\in M_n^{Lc}} p_{kn}\bigl(\delta_{k0}(n)-E_n
\bigl(\delta_{k0}(n)\bigr)\bigr) \biggr) \biggr)
\\[-4pt]
&&\qquad\leq\sum_{k\in M_n^{Lc}} \biggl( \log \biggl( \biggl(\exp
\biggl\{ \frac{r
na(b(n))}{b(n)} p_{kn} \biggr\}-1 \biggr) (1-
p_{kn})^n+1 \biggr)
\\[-4pt]
&&\hspace*{126.5pt}\qquad\quad{}-\frac{r na(b(n))}{ b(n)}p_{kn}(1- p_{kn})^n
\biggr).
\end{eqnarray*}
Therefore
\begin{eqnarray*}
&& \log E_n \biggl(\exp \biggl(\frac{r na(b(n))}{ b(n)}\sum
_{k\in M_n^{Lc}} p_{kn}\bigl(\delta_{k0}(n)-E_n
\bigl(\delta_{k0}(n)\bigr)\bigr) \biggr) \biggr)
\\[-4pt]
&&\qquad\leq 4\sum_{k\in M_n^{Lc}} \biggl(\frac{ r na(b(n))}{b(n)}
p_{kn} \biggr)^2e^{-n p_{kn}}I_{\{ {|r| na(b(n))}p_{kn}/{b(n)}\leq1 \}}
\\[-4pt]
&&\qquad\quad{} +12\sum_{k\in M_n^{Lc}} \exp \biggl\{\frac{2|r| na(b(n))}{b(n)}
p_{kn} \biggr\}e^{-n p_{kn}}I_{\{ {|r| na(b(n))}p_{kn}/{b(n)}\geq1
\}
}
\\[-4pt]
&&\qquad\leq \frac{4 r^2 a(b(n))}{b(n)} A_{nL}+ 24 |r| A_{n},\vadjust{\goodbreak}
\end{eqnarray*}
where $ A_{n}:=\frac{ a(b(n))}{b(n)}\sum_{k=1}^\infty p_{kn}
e^{-\lambda_n p_{kn}}I_{ \{\lambda_n p_{kn}\geq{b(n)}/({|r|
a(b(n))}) (1-{2|r| a(b(n))}/{b(n)}  )  \}}$,
$\lambda_n=n (1-\frac{2|r| a(b(n))}{b(n)}  )$ and $ A_{nL}:=
\frac{1}{b(n)} \sum_{k\in M_n^{Lc}} n^2 p_{kn}^2 e^{-n p_{kn}}. $ By
the proof of Lemma~\ref{truncation-lem-L}, $A_{n} \leq\frac{
a(b(n))}{b(n)}\exp \{-\frac{b(n)}{|r|a(b(n))} (1-\frac{2|r|
a(b(n))}{b(n)} )  \} \to0$. By\vspace*{2pt} (\ref{truncation-lem-L-eq-1}),
$ \limsup_{L\to\infty}\limsup_{n\to\infty}A_{nL}=0. $ Therefore,
(\ref{exp-moment-estimates-lem-1-eq-3}) holds.
\end{pf*}
\begin{pf*}{Proof of Lemma~\ref{laplace-int-lem}}
It is known that
\[
P_n \bigl(X_k(n) = x_k; k=1, \ldots,m \bigr) =P_n \Biggl(Y_k(n) = x_k;
k=1,\ldots,m \bigg|\sum_{k=1}^m Y_k(n)=n \Biggr),
\]
where $Y_k(n),k\geq1$ are independent random variables, and $Y_k(n)$
is Poisson distributed with mean $np_{kn}$.
Then, using the partial inversion formula for characteristic function
due to \citet{Bar} [see also \citet{Hol79}, \citet{Est83}], for any
$\alpha\in\mathbb R$,
\begin{eqnarray*}
&&
E_n \biggl(\exp \biggl\{\frac{\alpha a (b^L(n))}{b^L(n)} \sum
_{k\in
M_n^L}\bigl(\delta_{k1}(n)-np_{kn}
\delta_{k0}(n)\bigr) \biggr\} \biggr)
\\
&&\qquad=\frac{n!}{2\pi n^ne^{-n}}\\
&&\qquad\quad{}\times\int_{-\pi}^\pi
E_n \Biggl(\exp \Biggl\{ iu\sum_{l=1}^\infty
\bigl(Y_l(n)-np_{ln}\bigr)
\\
&&\hspace*{72.3pt}\qquad\quad{} +\frac{\alpha a (b^L(n))}{
b^L(n)}\sum_{k\in M_n^L}(I_{\{Y_k(n)=1\}}-np_{kn}I_{\{Y_k(n)=0\}
})
\Biggr\} \Biggr)\,du
\\
&&\qquad=\frac{n!}{2\pi n^ne^{-n}}\int_{-\pi}^{\pi}
H_{n}(u,\alpha)\,du,
\end{eqnarray*}
where
\begin{eqnarray*}
H_{n}(u,\alpha)&=&\prod_{k\in M_n^L} \bigl(
\theta_k(u,\alpha)+\gamma _k(u) \bigr) \prod
_{k\in M_n^{Lc}} \gamma_k(u),
\\
\gamma_k(u):\!&=&E_n \bigl(\exp \bigl\{iu
\bigl(Y_k(n)-np_{kn}\bigr) \bigr\} \bigr)= \exp \bigl
\{np_{kn}\bigl(e^{iu }-1-iu \bigr) \bigr\}
\end{eqnarray*}
and
\begin{eqnarray*}
\theta_k(u,\alpha) &=& \gamma_k(u)+\exp \{-iu n
p_{kn} \} \biggl(\exp \biggl\{ -\frac
{\alpha n a(b^L(n))}{b^L(n)} p_{kn}
\biggr\}-1 \biggr)\exp\{ -np_{kn}\}
\\
&&{}+\exp \bigl\{iu(1- n p_{kn})\bigr\}  \biggl(\exp \biggl\{\frac{\alpha a
(b^L(n))}{b^L(n)}
\biggr\}-1 \biggr)np_{kn}\exp\{-np_{kn}\}.
\end{eqnarray*}

It is obvious that $H_n(-u,\alpha)=\overline{H_n(u,\alpha)}$.
By Stirling's formula,
\[
\lim_{n\to\infty} \frac{n^ne^{-n}\sqrt {n}}{n!}=\frac{1}{\sqrt{2\pi}},
\]
it suffices to show that for any
$\alpha\in\mathbb R$,
%
\begin{equation}
\label{Laplace-int-lem-eq-2} \lim_{n\to\infty} \frac{ b^L(n)}{ a^2 (b^L(n))}\log
\int_{-\pi
}^{\pi
}n ^{1/2} H_n(u,
\alpha) \,du=\frac{\alpha^2}{2}.
\end{equation}

Since $ n p_{kn}\leq L $ uniformly in $k\in M_n^L$, we can write that
for $n$ large enough,
\[
H_n(u,\alpha) = \prod_{k\ge1}
\gamma_k(u)\prod_{k\in M_n^L} \bigl(1+
\gamma_k(u)^{-1}\theta_k(u,\alpha) \bigr)
=e^{n(e^{i u}-1-iu)}\prod_{k\in M_n^L} h_k(u,
\alpha),
\]
where
$ h_k(u,\alpha):=1+\gamma_k(u)^{-1}\theta_k(u,\alpha)$.

Choose a positive function $\kappa(t)$ such that $\kappa(t)\to\infty$
and $a(t)\kappa(t)/t\to0$, and define $\tau(t)= \sqrt{\frac
{a(t)(\kappa
(t))^{1/2}}{t}}$, $t\geq1 $ and then
$
\lim_{t\to\infty} \tau(t)=0$,\break $ \lim_{t\to\infty} \frac{\tau
^2(t)t}{a(t)}=\infty.
$
Noting that for $n$ large enough, $\sup_{u\in[\tau(n),\pi]}(1-\cos
u)\geq\tau^2(n)/4$, we have that
\begin{eqnarray*}
&& \frac{ b^L(n)}{ a^2 (b^L(n))}\log \Bigl( n ^{1/2}\sup_{ u \in
(\tau
(n),\pi]}\bigl|H_n(u,
\alpha)\bigr| \Bigr)
\\
&&\qquad\leq \frac{ b^L(n)\log n}{ 2 a^2 (b^L(n))}- \frac{ b^L(n)n \tau
^2(n)}{ 4 a^2 (b^L(n))}+\frac{ b^L(n)}{ a^2 (b^L(n))}\sum
_{k\in
M_n^L}\log\sup_{ u \in(\tau(n),\pi]}\bigl|h_k(u,
\alpha)\bigr|
\\
&&\qquad= - \frac{ b^L(n)n \tau^2(n)}{a^2(b^L(n))} \biggl(1+O \biggl(\frac
{\log
n}{n\tau^2(n)} \biggr)+O
\biggl(\frac{a(n)}{n\tau^2(n)} \biggr) \biggr) \to -\infty,
\end{eqnarray*}
which implies that
%
\begin{equation}
\label{Laplace-int-lem-eq-3} \limsup_{n\to\infty} \frac{ b^L(n)}{ a^2 (b^L(n))}
\log\biggl\llvert \int_{|u|\in[\tau(n),\pi]} n ^{1/2}H_n(u,
\alpha) \,du\biggr\rrvert =-\infty.
\end{equation}
Therefore, it suffices to show that
%
\begin{equation}
\label{Laplace-int-lem-eq-4} \limsup_{n\to\infty} \frac{ b^L(n)}{ a^2 (b^L(n))}
\log\int_{-\tau
(n)}^{\tau(n)} n ^{1/2}H_n(u,
\alpha) \,du=\frac{\alpha^2}{2}.
\end{equation}

In order to show (\ref{Laplace-int-lem-eq-4}), let us define a
transformation as follows. For $\alpha\in\mathbb R$ given, set $\rho
(n)= \frac{\alpha a(b^L(n))}{b^L(n) }\frac{E_n(F_1^L(n))}{n}$, and define
\[
\tilde{H}_n(z)=H_n \bigl(z+i\rho(n),\alpha \bigr),\qquad z
\in \mathbb C,
\]
where $\mathbb C$ denotes the complex plane. The transformation plays
an important role.
Let $\Gamma$ denote the closed path formed by the ordered points
$-\tau
(n)-i\rho(n)$, $\tau(n)-i\rho(n)$, $\tau(n)$, $-\tau(n)$, $-\tau
(n)-i\rho(n)$ on the complex plane.\vadjust{\goodbreak} Then by Cauchy's formula,
\begin{eqnarray*}
\int_{-\tau(n)}^{\tau(n)} H_n(u,\alpha) \,du &=&
\int_{-\tau(n)-i\rho(n)}^{\tau(n)-i\rho(n)} \tilde{H}_n(z) \,dz
\\
&=&-\int_{\tau(n)}^{-\tau(n)} \tilde{H}_n(z) \,dz-
\int_{\tau
(n)-i\rho
(n)}^{\tau(n)} \tilde{H}_n(z)\,dz\\
&&{}-\int
_{-\tau(n)}^{-\tau(n)-i\rho(n)} \tilde{H}_n(z) \,dz.
\end{eqnarray*}
Noting that
$\llvert \int_{\tau(n)-i\rho(n)}^{\tau(n)} \tilde{H}_n(z)\,dz\rrvert \leq
\llvert  \int_0^{\rho(n)}\tilde{H}_n (\tau(n)-iu
)\,du\rrvert $,
by
\begin{eqnarray*}
&&
\sup_{|u|\leq\rho(n)}\bigl\llvert \exp\bigl\{n\bigl(e^{-u}e^{i\tau(n)}-1-i
\tau (n)+ u\bigr)\bigr\} \bigr\rrvert \\
&&\qquad\leq\exp \biggl\{-n \biggl(
\frac{\tau^2(n)}{4}\bigl(1-\bigl|\rho(n)\bigr|\bigr)- \rho ^2(n) \biggr) \biggr\}
\end{eqnarray*}
and $\sup_{ |u|\leq|\rho(n)|}|h_k(\tau(n)+iu,\alpha)|=1+O (
\frac
{a(n)}{n}  ) n p_{kn}e^{-np_{kn}}$,
similarly to the proof of (\ref{Laplace-int-lem-eq-3}), we have that
\[
\frac{ b^L(n)}{ a^2 (b^L(n))}\log \biggl( n ^{1/2}\biggl\llvert \int
_{\tau
(n)-i\rho(n)}^{\tau(n)} \tilde{H}_n(z) \,dz\biggr
\rrvert \biggr) \to -\infty.
\]
Similarly,
$\frac{ b^L(n)}{ a^2 (b^L(n))}\log ( n ^{1/2}\llvert \int_{-\tau
(n)}^{-\tau(n)-i\rho(n)}\tilde{H}_n(z)\,dz\rrvert  ) \to-\infty$.
Therefore, it suffices to prove that
%
\begin{equation}
\label{Laplace-int-lem-eq-5} \limsup_{n\to\infty} \frac{ b^L(n)}{ a^2 (b^L(n))}
\log\int_{-\tau
(n)}^{\tau(n)} n ^{1/2}
\tilde{H}_n(u) \,du =\frac{\alpha^2}{2}.
\end{equation}
Let $\mathfrak{Re}(z)$ and $\mathfrak{Im}(z)$ denote the real part and
the imaginary part of a complex number $z$, respectively. Then
\begin{eqnarray*}
&&\mathfrak{Re} \bigl(h_k \bigl(u+i\rho(n),\alpha \bigr) \bigr)
\\
&&\qquad= 1+e^{np_{kn}(1-e^{-\rho(n)}\cos u)} \biggl(\cos \bigl(np_{kn}e^{-\rho
(n)}\sin u
\bigr) \\
&&\qquad\quad\hspace*{102pt}{}\times\biggl(\exp \biggl\{ - \frac{\alpha
a(b^L(n))}{b^L(n)} n p_{kn} \biggr\}-1
\biggr)e^{-np_{kn}}
\\
&&\hspace*{104pt}\qquad\quad{}+ \bigl(\cos \bigl(np_{kn}e^{-\rho(n)}\sin u
\bigr)e^{-\rho
(n)}\cos u\\
&&\hspace*{119pt}\qquad\quad{}+ \sin \bigl(np_{kn}e^{-\rho(n)}\sin u
\bigr)e^{-\rho(n)}\sin u \bigr)
\\
&&\hspace*{103.5pt}\hspace*{11.6pt}\qquad\quad{}\times \biggl(\exp \biggl\{ \frac{\alpha a(b^L(n))}{b^L(n)} \biggr\} -1
\biggr)np_{kn}e^{-np_{kn}} \biggr)
\end{eqnarray*}
and
\begin{eqnarray*}
&&\mathfrak{Im} \bigl(h_k \bigl(u+i\rho(n),\alpha \bigr) \bigr)
\\
&&\qquad= e^{np_{kn} (1-e^{-\rho(n)}\cos u )} \biggl(-\sin \bigl(np_{kn}e^{-\rho(n)}\sin u
\bigr) \\
&&\qquad\quad\hspace*{93.7pt}{}\times\biggl(\exp \biggl\{-\frac{\alpha
a(b^L(n))}{b^L(n)}n p_{kn} \biggr\}-1
\biggr)e^{-np_{kn}}
\\
&&\hspace*{84.3pt}\qquad\quad{} + \bigl(-\sin \bigl(np_{kn}e^{-\rho(n)}\sin u
\bigr)e^{-\rho
(n)}\cos u\\
&&\hspace*{84.3pt}\qquad\quad\hspace*{14.6pt}{}+\cos \bigl(np_{kn}e^{-\rho(n)}\sin u
\bigr)e^{-\rho(n)}\sin u \bigr)
\\
&&\hspace*{84.3pt}\hspace*{18pt}\qquad\quad{}\times \biggl(\exp \biggl\{ \frac{\alpha a(b^L(n))}{b^L(n)} \biggr\} -1
\biggr)np_{kn}e^{-np_{kn}} \biggr).
\end{eqnarray*}

For convenience, let $O_{jn}(u)$, $j\geq1$, denote uniformly bounded
real functions such that $O_{jn}(u)=0$ for all $|u|>\tau(n)$, and
$\lim_{n\to\infty}\sup_{u\in\mathbb R}|O_{jn}(u)|=0$.
Then for $n$ large enough, for all $u\in[-\tau(n),\tau(n)]$,
\begin{eqnarray*}
&&
\mathfrak{Re} \bigl(h_k \bigl(u+i\rho(n),\alpha \bigr) \bigr)
\\
&&\qquad= 1+e^{np_{kn}(1-e^{-\rho(n)}\cos u)} \biggl( \frac{1}{2} \biggl(\frac
{\alpha a(b^L(n))}{b^L(n)}
\biggr)^2\\
&&\qquad\quad\hspace*{102pt}{}\times \biggl(n^2 p_{kn}^2+
\biggl(1-\frac{2
E_n(F_1^L(n))}{n} \biggr)np_{kn} \biggr)e^{-np_{kn}}
\\
&&\hspace*{102pt}\qquad\quad{} +o \biggl(\frac{ a(n)}n \biggr)^2 n p_{kn}e^{-np_{kn}}
\\
&&\hspace*{157.6pt}\qquad\quad{}+ u^2 O_{1n}(u)O \biggl( \frac{ a (n)}{n}
\biggr)np_{kn}e^{-np_{kn}} \biggr)
\end{eqnarray*}
and
\begin{eqnarray*}
&&
\mathfrak{Im} \bigl(h_k \bigl(u+i\rho(n),\alpha \bigr) \bigr) \\
&&\qquad=
\frac{\alpha a(b^L(n))}{b^L(n)} e^{np_{kn}(1-e^{-\rho(n)}\cos
u)}un p_{kn}e^{-np_{kn}} \bigl(1+
u^2 O_{2n}(u) \bigr).
\end{eqnarray*}
Therefore
\begin{eqnarray*}
&&\bigl|H_n \bigl(u+i\rho(n),\alpha \bigr)\bigr|\\
&&\qquad = e^{-n(1-e^{-\rho(n)}\cos u-\rho(n))}\exp
\biggl\{ \frac{1}{2} \sum_{k\in M_n^L}
\log\bigl|h_k\bigl(u+i\rho(n),\alpha\bigr)\bigr|^2 \biggr\}
\\
&&\qquad=\exp \biggl\{\frac{1}{2}\frac{\alpha^2 a^2(b^L(n))}{b^L(n)}+o \biggl(
\frac
{ a^2(n)}{n} \biggr) \biggr\}\exp \biggl\{-\frac{n}{2}u^2
\bigl(1+ O_{4n}(u)o(1) \bigr) \biggr\}
\end{eqnarray*}
and so
\begin{eqnarray*}
&&
\frac{ b^L(n)}{ a^2 (b^L(n))}\log\int_{\tau(n)}^{-\tau(n)} n
^{1/2} \tilde{H}_n(u) \,du
\\
&&\qquad=\frac{\alpha^2}{2} +\frac{ b^L(n)}{ a^2 (b^L(n))} \log\int_{-\tau
(n)\sqrt{n}}^{\tau(n)\sqrt{n}}
\exp \biggl\{-\frac{1}{2}\frac
{\alpha^2
a^2(b^L(n))}{b^L(n)} \biggr\}\tilde{H}_n
\bigl(un^{-1/2} \bigr)\,du.
\end{eqnarray*}
Now, by
\begin{eqnarray*}
&&\frac{ b^L(n)}{ a^2 (b^L(n))} \log\int_{-\tau(n)\sqrt{n}}^{\tau
(n)\sqrt{n}}\exp \biggl
\{-\frac{1}{2}\frac{\alpha^2 a^2(b^L(n))}{b^L(n)} \biggr\}\bigl\llvert
\tilde{H}_n \bigl(un^{-1/2} \bigr)\bigr\rrvert \,du
\\
&&\qquad=o (1 )+\frac{ b^L(n)}{ a^2 (b^L(n))} \log\int_{-\tau
(n)\sqrt
{n}}^{\tau(n)\sqrt{n}}
\exp \biggl\{-\frac{1}{2}u^2 \bigl(1+ O_{4n}
\bigl(un^{-1/2}\bigr)o(1) \bigr) \biggr\}\,du\\
&&\qquad\to0,
\end{eqnarray*}
we obtain (\ref{Laplace-int-lem-eq-5}).
The proof of Lemma~\ref{laplace-int-lem} is complete.
\end{pf*}
\begin{pf*}{Proof of Lemma~\ref{exp-eqv-lemma-1}}
For any $t>\lambda
_n$, we can write [cf. (A.1) in \citet{ZhaZha09}]
%
\begin{eqnarray}
\label{exp-eqv-lemma-1-eq-2}
&&
Y_{ktn}-Y_{k\lambda_nn}\nonumber\\
&&\qquad=-Y_{k\lambda_n
n}I_{\{X_k(N_{t})>X_k(N_{\lambda_n})\}}\\
&&\qquad\quad{}+\delta_{k0}(N_{\lambda_n})
\bigl(\delta_{k1}(N_t)-(t-\lambda_n)p_{kn}
\delta_{k0}(N_t) \bigr).\nonumber
\end{eqnarray}
Therefore, it suffices to prove that
%
\begin{eqnarray}
\label{exp-eqv-lemma-1-eq-3}
&&\limsup_{n\to\infty}\frac{b(n)}{a^2(b(n))}\nonumber\\
&&\quad\hspace*{0pt}{}\times
\log P_n \Biggl(\sup_{\lambda_n\leq t\leq\lambda_n +\Delta_n} \sum
_{
k=1}^\infty Y_{k\lambda_n
n} I_{\{X_k(N_{t})>X_k(N_{\lambda_n})\}}\geq
\varepsilon a\bigl(b(n)\bigr) \Biggr)\\
&&\qquad=-\infty\nonumber
\end{eqnarray}
and
%
\begin{eqnarray}
\label{exp-eqv-lemma-1-eq-4}
&&\limsup_{n\to\infty}\frac{b(n)}{a^2(b(n))}\nonumber\\
&&\quad{}\times
\log P_n \Biggl(\sup_{\lambda
_n\leq t\leq\lambda_n +\Delta_n}\Biggl|\sum
_{ k=1}^\infty\delta _{k0}(N_{\lambda_n})
\bigl(\delta_{k1}(N_t)
-(t-\lambda_n)p_{kn}\delta_{k0}(N_t)
\bigr) \Biggr|\nonumber\\[-8pt]\\[-8pt]
&&\qquad\hspace*{242pt}\geq \varepsilon a\bigl(b(n)\bigr) \Biggr)\nonumber\\
&&\qquad=-\infty.
\nonumber
\end{eqnarray}

Let us first prove (\ref{exp-eqv-lemma-1-eq-3}). Set
$T_k=\min\{t\geq0;X_k(N_t)>X_k(N_{\lambda_n})\}$ and
$Z_t^{(n)}=\sum_{T_k\leq t}Y_{k\lambda_n n}$. Since $Y_{k\lambda_n
n},k\geq1$ are independent variables with mean zero and independent
of $\mathcal G:=\sigma({\mathbf X}(N_t)-{\mathbf
X}(N_{\lambda_n}),t\geq\lambda_n)$, $\{Z_t^{(n)},t\geq
\lambda_n\}$ is a martingale, and by the maximal inequality for
supermartingales, we have that for any $\varepsilon>0$, for any
$r>0$,
\begin{eqnarray*}
&&
P_n \Bigl(\sup_{\lambda_n\leq t\leq\lambda_n
+\Delta_n}\bigl|Z_t^{(n)}\bigr|
\geq\varepsilon a\bigl(b(n)\bigr) \Bigr)
\\
&&\qquad\leq 2e^{-r
\varepsilon a(b(n))}\max \bigl\{E_n \bigl(\exp \bigl\{r
Z_{\lambda_n
+\Delta_n}^{(n)} \bigr\} \bigr),E_n \bigl(\exp \bigl
\{-r Z_{\lambda_n
+\Delta_n}^{(n)} \bigr\} \bigr) \bigr\}
\end{eqnarray*}
and
\begin{eqnarray*}
&&
E_n \bigl(\exp \bigl\{r Z_{\lambda_n +\Delta_n}^{(n)} \bigr\}
\bigr) \\
&&\qquad= E_n \Biggl(E_n \Biggl(\exp \Biggl\{r \sum
_{k=1}^\infty Y_{k\lambda_n
n}I_{\{
X_k(N_{\lambda_n+\Delta_n})>X_k(N_{\lambda_n})\}}
\Biggr\} \bigg|\mathcal G \Biggr) \Biggr)
\\
&&\qquad=\prod_{k=1}^\infty \bigl( \bigl(
\bigl(e^{-r \lambda
_np_{kn}}+\lambda _np_{kn}e^{r}-1-
\lambda_np_{kn} \bigr)e^{-\lambda
_np_{kn}}\bigl(1-e^{-\Delta
_n p_{kn}}
\bigr)+ 1 \bigr) \bigr).
\end{eqnarray*}

For any $\alpha\not=0$, take $r= \frac{\alpha a(b(n))}{b(n)}$. Then for
$n$ large enough,
\begin{eqnarray*}
&& \exp \biggl\{- \frac{\alpha a(b(n))}{b(n)} \lambda_np_{kn}
\biggr\} +\lambda_np_{kn}\exp \biggl\{ \frac{\alpha a(b(n))}{b(n)}
\biggr\} -1-\lambda_np_{kn}
\\
&&\qquad\leq \frac{3\alpha^2 a^2(b(n))}{b^2(n)} \lambda_n^2p_{kn}^2I_{\{
{|\alpha| a(b(n))} \lambda_np_{kn}/{b(n)}\leq1\}}\\
&&\qquad\quad{}
+ \exp \biggl\{ \frac{
|\alpha| a(b(n))}{b(n)}\lambda_n p_{kn} \biggr
\}I_{\{
{|\alpha| a(b(n))}\lambda_np_{kn}/{b(n)}> 1\}}
\\
&&\qquad\quad{}+ \frac{3\alpha^2 a^2(b(n))}{b^2(n)} \lambda_np_{kn}.
\end{eqnarray*}
Therefore, for $n$ large enough,
\begin{eqnarray*}
&& \sum_{k=1}^\infty \biggl(\exp \biggl\{-
\frac{\alpha a(b(n))}{b(n)} \lambda_np_{kn} \biggr\}+
\lambda_np_{kn}\exp \biggl\{ \frac{\alpha
a(b(n))}{b(n)} \biggr
\}-1-\lambda_np_{kn} \biggr)\\
&&\quad\hspace*{2pt}{}\times e^{-\lambda_np_{kn}
}
\bigl(1-e^{-\Delta_n p_{kn}}\bigr)
\\
&&\hspace*{2pt}\qquad\leq \frac{3\alpha^2 a^2(b(n))}{b (n)} B_{1n}+ \frac{\alpha^2 a^2
(b(n))}{b(n)}
B_{2n} +\frac{3\alpha^2 a^2(b(n))}{b(n)}B_{3n},
\end{eqnarray*}
where
$
B_{1n}:=\frac{1}{b(n)}\sum_{k=1}^\infty\lambda_n^2p_{kn}^2I_{\{
{|\alpha| a(b(n))}\lambda_np_{kn}/{b(n)}\leq1\}} e^{-\lambda_np_{kn}
}(1-e^{-\Delta_n p_{kn}}),
$
\[
\hspace*{-10pt}B_{2n}:= \frac{1}{b(n)}\sum_{k=1}^\infty
\sum_{k=1}^\infty\lambda_n^2p_{kn}^2
\exp \biggl\{- \biggl(1-\frac{|\alpha| a(b(n))}{b(n)} \biggr) \lambda_n
p_{kn} \biggr\}I_{\{
\lambda_np_{kn}> {b(n)}/({|\alpha| a(b(n))})\}}
\]
and
$
B_{3n}:=\frac{1}{b(n)}\sum_{k=1}^\infty\lambda_np_{kn} e^{-\lambda
_np_{kn} }(1-e^{-\Delta_n p_{kn}}).
$
By (\ref{mdp-Lindeberg-condition}),
$
B_{2n}\to0.
$
Then, by $\frac{s_n^2}{n}\to0$ under $c_1=0$, $ s_{n}^2/b(n)\to1$ and
$ s_{\lambda_nn}^2/s_n^2\to1$,
\begin{eqnarray*}
B_{1n} &\leq&\frac{4M}{|\alpha|}\sqrt{\frac{1}{\lambda_nb(n)}}\sum
_{k=1}^\infty\lambda_n^2p_{kn}^2I_{\{ {|\alpha|a(b(n))}
\lambda_np_{kn}/{b(n)}\leq1\}}
e^{-\lambda_np_{kn} } \\
&\leq&\frac{8 M}{|\alpha| } s_{\lambda_nn}^2\sqrt{
\frac{1}{\lambda
_nb(n)}}\to0
\end{eqnarray*}
and
\begin{eqnarray*}
B_{3n}&\leq&\frac{4M}{|\alpha| }\sqrt{\frac{1}{\lambda_nb(n)}}\sum
_{k=1}^\infty\lambda_np_{kn}I_{\{ {|\alpha| a(b(n))}
\lambda_np_{kn}/{b(n)}\leq1\}}
e^{-\lambda_np_{kn} }
\\
&&{}+\frac{2}{b(n)}\sum_{k=1}^\infty
\lambda_np_{kn}e^{-\lambda
_np_{kn}}I_{\{{|\alpha| a(b(n))}\lambda_np_{kn}/{b(n)}>1\}}\to0.
\end{eqnarray*}
Thus
%
\begin{equation}
\label{exp-eqv-lemma-1-eq-5} \lim_{n\to\infty}\frac{b(n)}{a^2(b(n))}\log
E_n \biggl(\exp \biggl\{\frac{\alpha a(b(n))}{b(n)} Z_{\lambda_n
+\Delta_n}^{(n)}
\biggr\} \biggr)=0.
\end{equation}
This yields (\ref{exp-eqv-lemma-1-eq-3}) by Chebyshev's inequality.
Next, we show (\ref{exp-eqv-lemma-1-eq-4}). Noting that
\begin{eqnarray*}
&&\sup_{\lambda_n\leq t\leq\lambda_n
+\Delta_n}\bigl|\delta_{k0}(N_{\lambda_n}) \bigl(
\delta_{k1}(N_t)-(t-\lambda_n)p_{kn}
\delta_{k0}(N_t) \bigr)\bigr|
\\
&&\qquad\leq \delta_{k0}(N_{\lambda_n}) (I_{\{X_k(N_{\lambda_n
+\Delta_n})>X_k(N_{\lambda_n})\}} +
\Delta_n p_{kn} ),
\end{eqnarray*}
it suffices to show that for any $\varepsilon>0$,
%
\begin{eqnarray}
\label{exp-eqv-lemma-1-eq-6}\qquad
&&\lim_{n\to\infty}\frac{b( n)}{a^2(b(n))}\log
P_n \Biggl(\sum_{k=1}^\infty
\delta_{k0}(N_{\lambda_n})I_{\{
X_k(N_{\lambda_n
+\Delta_n})>X_k(N_{\lambda_n})\}} >\varepsilon a\bigl(b(n)
\bigr) \Biggr)\nonumber\\[-8pt]\\[-8pt]
&&\qquad=-\infty\nonumber
\end{eqnarray}
and
%
\begin{equation}
\label{exp-eqv-lemma-1-eq-7} \lim_{n\to\infty}\frac{b( n)}{a^2(b(n))}\log
P_n \Biggl(\sqrt{\frac{n}{b(n)}}\sum
_{k=1}^\infty\delta _{k0}(N_{\lambda
_n})
p_{kn} >\varepsilon \Biggr)=-\infty.
\end{equation}
Since
\begin{eqnarray*}
&&E_n \Biggl(\exp \Biggl\{ \frac{\alpha a (b(n))}{b(n)}\sum
_{k=1}^\infty \delta_{k0}(N_{\lambda_n})I_{\{X_k(N_{\lambda_n +\Delta
_n})>X_k(N_{\lambda_n})\}}
\Biggr\} \Biggr)
\\
&&\qquad= \prod_{k=1}^\infty \biggl( 1 + \biggl(\exp
\biggl\{ \frac{\alpha a
(b(n))}{b(n)} \biggr\}-1 \biggr) \bigl(1-e^{-\Delta_n p_{kn}}
\bigr)e^{-\lambda_n
p_{kn}} \biggr),
\end{eqnarray*}
a similar argument to the proof of (\ref{exp-eqv-lemma-1-eq-5}) gives
\[
\hspace*{-6pt}\lim_{n\to\infty}\frac{b(n)}{ a ^2(b(n))}\log E_n \!\Biggl(
\exp \!\Biggl\{ \frac{\alpha a
(b(n))}{b(n)}\sum_{k=1}^\infty\!
\delta_{k0}(N_{\lambda_n})I_{\{
X_k(N_{\lambda_n
+\Delta_n})>X_k(N_{\lambda_n})\}}\! \Biggr\}\! \Biggr)=0,
\]
which implies (\ref{exp-eqv-lemma-1-eq-6}).
Similarly, we can obtain (\ref{exp-eqv-lemma-1-eq-7}).
\end{pf*}
\end{appendix}

\section*{Acknowledgments}

The author is very grateful to the Editor, Professor T. Cai, the
Associate Editor and an anonymous referee for their helpful comments
and suggestions. The author is also thankful to the referee for
recommending the reference \citet{LlaGouRee11}.



\printaddresses

\end{document}